\documentclass{article}
\usepackage{amssymb,amsmath,theorem,euscript}

\newcounter{sec}

\def\sm{\smallskip}

\newcounter{punct}[sec]

\def\punct{\refstepcounter{punct}{\arabic{sec}.\arabic{punct}.  }}

\def\COUNTERS{\addtocounter{sec}{1}
              \setcounter{punct}{0}
          \setcounter{equation}{0}
          \setcounter{theorem}{0}
                  }

\newtheorem{theorem}{Theorem}[sec]
\newtheorem{proposition}[theorem]{Proposition}
\newtheorem{lemma}[theorem]{Lemma}
\newtheorem{corollary}[theorem]{Corollary}

\begin{document}

 \def\ov{\overline}
\def\wt{\widetilde}
 \newcommand{\rk}{\mathop {\mathrm {rk}}\nolimits}
\newcommand{\Aut}{\mathop {\mathrm {Aut}}\nolimits}
\newcommand{\Out}{\mathop {\mathrm {Out}}\nolimits}
\renewcommand{\Re}{\mathop {\mathrm {Re}}\nolimits}
\def\Br{\mathrm {Br}}

 %%%1. ClASSICAL GROUPS
\def\SL{\mathrm {SL}}
\def\SU{\mathrm {SU}}
\def\GL{\mathrm {GL}}
\def\U{\mathrm U}
\def\OO{\mathrm O}
 \def\Sp{\mathrm {Sp}}
 \def\SO{\mathrm {SO}}
\def\SOS{\mathrm {SO}^*}
 \def\Diff{\mathrm{Diff}}
 \def\Vect{\mathfrak{Vect}}
\def\PGL{\mathrm {PGL}}
\def\PU{\mathrm {PU}}
\def\PSL{\mathrm {PSL}}
\def\Symp{\mathrm{Symp}}
\def\End{\mathrm{End}}
\def\Mor{\mathrm{Mor}}
\def\Aut{\mathrm{Aut}}
 \def\PB{\mathrm{PB}}
 \def\cA{\mathcal A}
\def\cB{\mathcal B}
\def\cC{\mathcal C}
\def\cD{\mathcal D}
\def\cE{\mathcal E}
\def\cF{\mathcal F}
\def\cG{\mathcal G}
\def\cH{\mathcal H}
\def\cJ{\mathcal J}
\def\cI{\mathcal I}
\def\cK{\mathcal K}
 \def\cL{\mathcal L}
\def\cM{\mathcal M}
\def\cN{\mathcal N}
 \def\cO{\mathcal O}
\def\cP{\mathcal P}
\def\cQ{\mathcal Q}
\def\cR{\mathcal R}
\def\cS{\mathcal S}
\def\cT{\mathcal T}
\def\cU{\mathcal U}
\def\cV{\mathcal V}
 \def\cW{\mathcal W}
\def\cX{\mathcal X}
 \def\cY{\mathcal Y}
 \def\cZ{\mathcal Z}
%%% END MATHCAL %%%%%%%%%%%%%%%%%%%%%%%%%%%%%%%%% %%%%%%%%%%%%%%%%%%%%%%%%%%%%%%%% %%%
\def\0{{\ov 0}}
 \def\1{{\ov 1}}
 %%%%%%%%%%%%%%%%%%%%%%%%%%%% %%%%%%%%%%%%%%%%%%%%%%%%%%%%%%%%%%% %%% BEGIN GOTIC
 \def\frA{\mathfrak A}
 \def\frB{\mathfrak B}
\def\frC{\mathfrak C}
\def\frD{\mathfrak D}
\def\frE{\mathfrak E}
\def\frF{\mathfrak F}
\def\frG{\mathfrak G}
\def\frH{\mathfrak H}
\def\frI{\mathfrak I}
 \def\frJ{\mathfrak J}
 \def\frK{\mathfrak K}
 \def\frL{\mathfrak L}
\def\frM{\mathfrak M}
 \def\frN{\mathfrak N} \def\frO{\mathfrak O} \def\frP{\mathfrak P} \def\frQ{\mathfrak Q} \def\frR{\mathfrak R}
 \def\frS{\mathfrak S} \def\frT{\mathfrak T} \def\frU{\mathfrak U} \def\frV{\mathfrak V} \def\frW{\mathfrak W}
 \def\frX{\mathfrak X} \def\frY{\mathfrak Y} \def\frZ{\mathfrak Z} \def\fra{\mathfrak a} \def\frb{\mathfrak b}
 \def\frc{\mathfrak c} \def\frd{\mathfrak d} \def\fre{\mathfrak e} \def\frf{\mathfrak f} \def\frg{\mathfrak g}
 \def\frh{\mathfrak h} \def\fri{\mathfrak i} \def\frj{\mathfrak j} \def\frk{\mathfrak k} \def\frl{\mathfrak l}
 \def\frm{\mathfrak m} \def\frn{\mathfrak n} \def\fro{\mathfrak o} \def\frp{\mathfrak p} \def\frq{\mathfrak q}
 \def\frr{\mathfrak r} \def\frs{\mathfrak s} \def\frt{\mathfrak t} \def\fru{\mathfrak u} \def\frv{\mathfrak v}
 \def\frw{\mathfrak w} \def\frx{\mathfrak x} \def\fry{\mathfrak y} \def\frz{\mathfrak z} \def\frsp{\mathfrak{sp}}
 %% This is Lie algebra %%% END GOTIC
%%%%%%%%%%%%%%%%%%%%%%%%%%%%%%%% %%%%%%%%%%%%%%%%%%%%%%%%%%%%%%%%%
%%% BEGIN MATHBF
 \def\bfa{\mathbf a} \def\bfb{\mathbf b} \def\bfc{\mathbf c} \def\bfd{\mathbf d} \def\bfe{\mathbf e} \def\bff{\mathbf f}
 \def\bfg{\mathbf g} \def\bfh{\mathbf h} \def\bfi{\mathbf i} \def\bfj{\mathbf j} \def\bfk{\mathbf k} \def\bfl{\mathbf l}
 \def\bfm{\mathbf m} \def\bfn{\mathbf n} \def\bfo{\mathbf o} \def\bfp{\mathbf p} \def\bfq{\mathbf q} \def\bfr{\mathbf r}
 \def\bfs{\mathbf s} \def\bft{\mathbf t} \def\bfu{\mathbf u} \def\bfv{\mathbf v} \def\bfw{\mathbf w} \def\bfx{\mathbf x}
 \def\bfy{\mathbf y} \def\bfz{\mathbf z} \def\bfA{\mathbf A} \def\bfB{\mathbf B} \def\bfC{\mathbf C} \def\bfD{\mathbf D}
 \def\bfE{\mathbf E} \def\bfF{\mathbf F} \def\bfG{\mathbf G} \def\bfH{\mathbf H} \def\bfI{\mathbf I} \def\bfJ{\mathbf J}
 \def\bfK{\mathbf K} \def\bfL{\mathbf L} \def\bfM{\mathbf M} \def\bfN{\mathbf N} \def\bfO{\mathbf O} \def\bfP{\mathbf P}
 \def\bfQ{\mathbf Q} \def\bfR{\mathbf R} \def\bfS{\mathbf S} \def\bfT{\mathbf T} \def\bfU{\mathbf U} \def\bfV{\mathbf V}
 \def\bfW{\mathbf W} \def\bfX{\mathbf X} \def\bfY{\mathbf Y} \def\bfZ{\mathbf Z} \def\bfw{\mathbf w}
 %%% END MATHBF
%%%%%%%%%%%%%%%%%%%%%%%%%%%%%%% %%%%%%%%%%%%%%%%%%%%%%%%%%%%%%%%%
 %%% BEGIN MATHBB
 \def\R {{\mathbb R }} \def\C {{\mathbb C }} \def\Z{{\mathbb Z}} \def\H{{\mathbb H}} \def\K{{\mathbb K}}
 \def\N{{\mathbb N}} \def\Q{{\mathbb Q}} \def\A{{\mathbb A}} \def\T{\mathbb T} \def\P{\mathbb P} \def\G{\mathbb G}
 \def\bbA{\mathbb A} \def\bbB{\mathbb B} \def\bbD{\mathbb D} \def\bbE{\mathbb E} \def\bbF{\mathbb F} \def\bbG{\mathbb G}
 \def\bbH{\mathbb H}
 \def\bbI{\mathbb I} \def\bbJ{\mathbb J} \def\bbL{\mathbb L} \def\bbM{\mathbb M} \def\bbN{\mathbb N} \def\bbO{\mathbb O}
 \def\bbP{\mathbb P} \def\bbQ{\mathbb Q} \def\bbS{\mathbb S} \def\bbT{\mathbb T} \def\bbU{\mathbb U} \def\bbV{\mathbb V}
 \def\bbW{\mathbb W} \def\bbX{\mathbb X} \def\bbY{\mathbb Y} \def\kappa{\varkappa} \def\epsilon{\varepsilon}
 \def\phi{\varphi} \def\le{\leqslant} \def\ge{\geqslant}

\def\UU{\bbU}
\def\Mat{\mathrm{Mat}}
\def\tto{\rightrightarrows}

\def\Gr{\mathrm{Gr}}

\def\graph{\mathrm{graph}}

\def\O{\mathrm{O}}

\def\la{\langle}
\def\ra{\rangle}

\renewcommand{\Im}{\mathop {\mathrm {Im}}\nolimits}
\renewcommand{\Re}{\mathop {\mathrm {Re}}\nolimits}

\def\sm{\smallskip}

\newcommand{\const}{\mathop {\mathrm {const}}\nolimits}
\newcommand{\res}{\mathop {\mathrm {res}}\nolimits}

\begin{center}
\Large\bf
Difference Sturm--Liouville problems in the imaginary direction

\bigskip

\large\sc
Yury A. Neretin%
\footnote{Supported by the grant FWF P22122.}

\end{center}

{\small
We consider difference operators in $L^2$ on $\R$  of the form
$$
\cL f(s)=p(s)f(s+i)+q(s) f(s)+r(s) f(s-i)
,$$
where $i$ is the imaginary unit. The domain of definiteness
are functions holomorphic in a strip with some conditions of decreasing
at infinity.
Problems of such type with discrete spectra are well known
(Meixner--Pollaszek, continuous Hahn, continuous dual Hahn,
and Wilson hypergeometric orthogonal polynomials).
 We write explicit spectral decompositions
for several operators $\cL$ with continuous spectra. 
We also discuss analogs of 'boundary conditions'
for such operators.}

%Key words Difference operators, Sturm--Liouville problem,
%Kontorovich--Lebedev transform, Wimp transform,
% index hypergeometric transform, defect indices, self-adjoint operator,
%spectral decomposition

%34B24    Sturm-Liouville theory 
%47B39     Difference operators 
%47G10      Integral operators
%44A20      Transforms of special functions 
%58C40  Spectral theory; eigenvalue problems 
%33C15   Confluent hypergeometric functions, Whittaker functions,
%33C05    Classical hypergeometric functions, 

\section{Introduction}

\COUNTERS

{\bf\punct Formulation of problem.}
Consider the space $L^2$ on $\R$ with respect to a positive weight $w(s)\,ds$. 
Consider a subspace
$H$ consisting of functions $f(s)$ holomorphic in the strip $-1<\Im s<1$ smooth up to the boundary
$\Im s=\pm 1$ and sufficiently rapidly decreasing in the strip as $|s|\to\infty$. 
We consider difference operators in $L^2(\R, w(s)\,ds)$  of the form
$$
\cL f(s)=p(s)f(s+i)+q(s) f(s)+r(s) f(s-i)
,$$
where $i$ is the imaginary unit; the domain of definiteness of $\cL$  is the subspace $H$.
For such operators we discuss essential self-adjointness and the eigenvalue
problem
$$
\cL f(s)=\lambda f(s)
.$$
Our main purpose is spectral decomposition.
In fact, several problems of this kind were solved (see the list below). All 
solved problems had the following form. Denote
\begin{equation}
\mu(s)=e^{cs}\frac{\prod_{k=1}^m \Gamma(a_k+ is)}{\prod_{l=1}^n\Gamma(b_l+is)}
\label{eq:mu}
,\end{equation}
where $c\in\R$, and
\begin{equation}
\nu(s)=\ov{\mu(\ov s)}=
e^{c s}\frac{\prod_{k=1}^m \Gamma(\ov a_k- is)}{\prod_{l=1}^n\Gamma(\ov b_l-is)}
\label{eq:nu}
.
\end{equation}
Denote
\begin{align}
A(s):=\frac{\nu(s+i)}{\nu(s)}= e^{-i c}\frac{\prod_{k=1}^m 
(\ov a_k- is)}{\prod_{l=1}^n(\ov b_l-is)}
,
\label{eq:A(s)}
\\
B(s):=\frac{\mu(s-i)}{\mu(s)}= e^{ic}\frac{\prod_{k=1}^m 
(a_k+ is)}{\prod_{l=1}^n(b_l+is)}
.
\label{eq:B(s)}
\end{align}
We consider the space $L^2(\R, w(s)\,ds)$ with respect to the weight
$$
w(s)\,ds:= \frac 1{2\pi}\mu(s)\nu(s)
$$
and the difference operator
\begin{equation}
\cL f(s)= A(s)f(s+i)-\bigl(A(s)+B(s)\bigr)f(s)+B(s)f(s-i)
\label{eq:difference-operator}
.
\end{equation}

%%%%%%%%%urm%%%%%%%%%%%%%%%%%%%%%%%%%%%%%%%%%%%%%%%%%%%%%%%%%%%%%%%%%%

{\bf\punct Neo-classical orthogonal polynomials}.
Now we enumerate solved problems of this kind.
We use the standard notation for hypergeometric functions
$$
_pF_q\left[\begin{matrix} a_1,\dots,a_q\\ b_1,\dots, b_p \end{matrix}\,;\,z\right]
:=\sum_{n=0}^\infty \frac{(a_1)_n \dots (a_p)_n\, z^n}{(b_1)_n\dots (b_q)_n\, n!}
,$$
where
$(a)_n:=a(a+1)\dots(a+n-1)$ is the Pochhammer symbol.

Recall that there are 3 types of classical hypergeometric orthogonal
 polynomials,
see \cite{AW}, \cite{Koe}, \cite{Koe2}.
Polynomials of the first  type are solutions
 of the usual Sturm--Liouville problems 
for second order differential operators: Jacobi (including Gegenbauer, Legendre, Chebyshev),
Laguerre, Hermite systems (see \cite{HTF2}). 

Polynomials of 
the second type are solutions of 
 difference Sturm--Liouville problemы on
lattices: Racah, (Chebyshev)--Hahn, dual Hahn, 
Meixner, Krawtchouk, Charlier, see \cite{NUS},  \cite{Koe}, \cite{Koe2}. 

Polynomials of the  third type
are solutions of Sturm--Liouville problems  of the form
(\ref{eq:mu})--(\ref{eq:difference-operator}): Wilson, continuous Hahn,
continuous dual Hahn, Meixner--Pollaczek systems,
 see \cite{Koe}, \cite{AAR}.
 Recall that all classical polynomial orthogonal
systems are degenerations of the Wilson polynomials,
see \cite{AW}, \cite{Koe}, \cite{Koe2}.

\sm

a) {\it The Meixner--Pollaczek system} or {\it Meixner
polynomials of the second kind}, see \cite{Mei}, \cite{Koe}, Section 1.7. We take
$$
\mu(s)=e^{(\phi-\pi/2)s}\Gamma(a+is)
,$$
where parameters $a$, $\phi$ satisfy $a> 0$, $0<\phi<\pi$.
Therefore
\begin{equation}
w(s)=\frac 1{2\pi}
e^{(2\phi-\pi)s}\Gamma(a+is)\Gamma(a-is)
.
\label{eq:measure-meixner}
\end{equation}
and the difference operator is
\begin{multline}
\cL f(s)
=
ie^{-i \phi}(a-is)f(s+i) + 2(-s\cos\phi+\lambda \sin\phi)f(s)-
\\=
ie^{i \phi}(a+is)f(s-i)
\label{eq:L-meixner}
.
\end{multline}

 The eigenfunctions are polynomials
\begin{align*}
P_n(s)&=\frac{(2a)_n}{n!} e^{in\phi}
 {}_2F_1\left[\begin{matrix} -n,a+is\\ 2a \end{matrix}; 1-e^{-2i\phi} \right],
\\ 
  \cL P_n(s)&= n\sin \phi\, P_n(s)
.
\end{align*}
Norms of Meixner--Pollaczek polynomials are given by
$$
\|w(s)\|^2:=\int_{-\infty}^{\infty} |p_n(s)|^2 w(s)\,ds=
\frac{\Gamma(n+2a)}{(2\sin\phi)\,n!}
.
$$

Recall (see \cite{HTF1}, formula 1.18(6))  that
\begin{equation}
|\Gamma(a+is)|\sim \sqrt{2\pi} |s|^{a-1/2} e^{-\pi s/2},\qquad s\to\infty
\label{eq:gamma-is}
.
\end{equation}
Therefore the weight $w(s)$  exponentially decreases and the space
$L^2(\R,w(s)\,ds)$ contains all polynomials. The operator $\cL$ send
 a polynomial to a polynomial of the same degree, therefore
 our Sturm--Liouville problem is pure algebraic. The same remarks 
 hold for 3 polynomial systems discussed below.

\sm

b) {\it The continuous Hahn system,} see \cite{Car}, \cite{AW0}, \cite{NUS}, \cite{Koe}.
In this case, 
$$
\mu(s)=\Gamma(a+is)\Gamma(b+is)
,$$
where the parameters $a$, $b$ satisfy $\Re a>0$, $\Re b>0$. The eigenfunctions 
are polynomials
\begin{align*}
p_n(s)&:=i^n \frac{(a+\ov a)_n(a+\ov b)_n}{n!}\,
{}_3F_2\left[\begin{matrix} -n, n+a+b+\ov a+\ov b,a+is\\
a+\ov a, a+\ov b \end{matrix} ; 1 \right],
\\
\cL p_n&=n(n+a+\ov a+b+\ov b) p_n
.
\end{align*}

\sm

c) {\it The continuous dual Hahn system,} see \cite{Wil}, \cite{Koe}, \cite{Koe2}.
In this case
$$
\mu(s)=\frac{\Gamma(a+is)\Gamma(b+is)\Gamma(c+is)}
{\Gamma(2is)}
,$$
where the parameters $a$, $b$, $c$ satisfy 
$a>0$, $b>0$, $c>0$ or $a>0$, $\Re b>0$, $c=\ov b$.
We consider {\it even} orthogonal polynomials $p_n(s^2)$: 
\begin{align*}
p_n(s^2)&:=(a+b)_n(a+c)_n{}\,\,{}_3F_2\left[\begin{matrix}-n,a+is,a-is\\ a+b,a+c \end{matrix};\,1  \right],
\\
\cL p_n&=np_n
.\end{align*}
%Norms of $p_n$ are given by
%$$
%\|p_n\|^2=\Gamma(n+a+b)\Gamma(n+a+c)\Gamma(n+b+c)
%.$$

{\it d) Wilson system,} see \cite{Wil}, \cite{AAR},
\cite{Koe}, \cite{Ner-wilson}.
   In this case,
$$
\mu(s)=\frac{\Gamma(a+is)\Gamma(b+is)\Gamma(c+is)\Gamma(d+is)}{\Gamma(2is)}
,$$
where
$\Re a$, $\Re b$, $\Re c$, $\Re d>0$
and all parameters are real, or $a$, $b$ are real, $d=\ov c$,
or $b=\ov a$, $d=\ov c$.
Wilson polynomials are {\it even} polynomials given by
\begin{multline*}
P_n(a,b,c,d;s^2)=\\=(a+b)_n (a+c)_n (a+d)_n
\,\,
{}_4F_3 \left[\begin{matrix}
-n,n+a+b+c+d-1, a+is,a-is\\
a+b,a+c,a+d\end{matrix}; 1\right]
.
\end{multline*}
They satisfy to the difference equation
$$
\cL P_n = n(a+b+c+d-1)P_n
.
$$
%Norms are 
%$$
%\|p_n\|^2=
%\frac{n!\Gamma(a+b+n)\Gamma(a+c+n)\Gamma(a+d+n)
%  \Gamma(b+c+n)\Gamma(b+d+n)\Gamma(c+d+n)}
%  {\Gamma(a+b+c+d+n) (a+b+c+d+2n-1)}    
%.$$

%%%%%%%%%%%%%%%%%%%%%%%%%%%%%%%%%%%%%%%%%%%%%%%%%%%%%%%%%%%%

{\bf\punct Sturm--Liouville problems with continuous spectra.}
I know two solved problems.

\sm

a) We consider {\it even} functions $f(s)$ on the line,
and
$$
\mu(s)=\frac{\Gamma(a+is)\Gamma(b+is)}{\Gamma(2is)}
,$$
where $a$, $b>0$. Let $\cL$ be the same as above.

We consider the operator (it is called the inverse Olevsky transform, \cite{Ole},
or the inverse Jacobi transform, \cite{Koo}):
$$
L^2\left(\R_+, \frac 1\pi
\left|\frac{\Gamma(a+is)\Gamma(b+is)}{\Gamma(2is)}\right|^2 ds\ \right)
\to
L^2(\R_+,x^{a+b-1}(1+x)^{a-b})
$$ 
defined by
$$
Jf(x)=\frac1{\pi\Gamma(a+b)}\int_0^\infty
{}_2F_1\left[\begin{matrix}a+is,a-is\\a+b \end{matrix}; x\right]
f(x)
\left|\frac{\Gamma(a+is)\Gamma(b+is)}{\Gamma(2is)}\right|^2 ds
.
$$
The $J$ send the difference operator $\cL$ to the operator
$$
Mf(x)=xf(x)
.
$$
See  \cite{Ner-ber}, Theorem 2.1, but this is very special case
of Cherednik, \cite{Che1}.

\sm

b) Let
$$
\mu(s)=\frac{\Gamma(a+is)\Gamma(b+is)\Gamma(c+is)}{\Gamma(d+is)\Gamma(2is)}
.
$$
In this case the spectral decomposition was done by an integral operator, whose
kernel is a $_4F_3$-function, see Groenevelt \cite{Gro},
the discrete part of the spectrum was  found in
\cite{Ner-wilson}.

\sm

%%%%%%%%%%%%%%%%%%%%%%%%%%%%%%%%%%%%%%%%%%%%%%%%%%%%%%%

{\bf\punct Partially solved problems. Romanovski-type systems of orthogonal polynomials.}
 Romanovski \cite{Rom} constructed orthogonal polynomials
on $\R$ with respect to the weight $(1+ix)^{-a}(1-ix)^{-\ov a}$ on $\R$
and  with respect the weight $x^{a-1}(1+x)^{-b}$ on $(0,\infty)$.
Since the weights have polynomial decreasing, these
orthogonal  systems are finite.
However, Romanovski polynomials correspond to discrete part
of spectra of certain  Sturm--Liouville problems
(see  \cite{Dunford}, XIII.8, \cite{Koo}, \cite{Ner-double}). 

Lesky (see, e.g., \cite{Les1}, \cite{Les2}) constructed numerous Romanovski type polynomial
systems related to difference Sturm--Liouville problems, his list contains
several difference problems in imaginary direction%
\footnote{More generally, Lesky's papers indicate numerous unsolved but
(certainly) solvable Sturm-Liouville
problems.}.

%%%%%%%%%%%%%%%%%%%%%%%%%%%%%%%%%%%%%%%%%%%%%%%%%%%%%%%%%%%%%%%

\sm

{\bf\punct Multidimensional analogs.} See \cite{Che1}, \cite{Che2}.

\sm

{\bf \punct Results of the paper.} In Section 2, we show that the operators
 (\ref{eq:mu})--(\ref{eq:difference-operator}) are formally symmetric.
 Next,
we found spectral decomposition
for several operators $\cL$.
In Sections 3, 4 we consider
$$
\mu(s)=\frac 1{\Gamma(is)}
\qquad\text{and}\qquad
\mu(s)=\frac{\Gamma(a+is)}{\Gamma(2is)}
$$
respectively.
In both cases the spectrum is the half-line $\lambda>0$.
The spectral decomposition is given respectively by the 
inverse Kontorovich--Lebedev
transform and the inverse Wimp transform with Whittaker kernel.
Note that in a certain sense these problems (involving Bessel functions
$_0F_1$ and Kummer functions $_1F_1$) are simpler than neo-classical
polynomial problems
(involving  the Gauss function $_2F_1$ and higher hypergeometric functions
 $_3F_2(1)$, $_4F_3(1)$).

Next (Section 5), we consider $L^2(\R)$ with respect to the measure
\begin{equation}
e^{\pi s}
| \Gamma(\alpha/2+is)|^2\,ds
\label{eq:measure-vil}
\end{equation}
and the difference operator
\begin{equation}
\cL f(s)=i(\alpha/2+is) f(s-i)+ 2 \cosh \phi \,s h(s)- i(\alpha/2-is) f(s+i)
.
\label{eq:L-vil}
\end{equation}
The form of this
 operator slightly differs from  (\ref{eq:mu})--(\ref{eq:difference-operator}).

 In Section 6 we discuss an example of a symmetric
non self-adjoint operator  and its essentially self-adjoint extensions.

\sm

In all cases essential self-adjointness is derived from the 
explicit spectral decomposition. It is an interesting question
to find a priory proofs.

\sm

We also note that the problem (\ref{eq:measure-vil})--(\ref{eq:L-vil})
is an analytic continuation of the Meixner--Pollaszek problem
 (\ref{eq:measure-meixner})--(\ref{eq:L-meixner}).
 The objects of Section 6 also are "analytic continuations
 from integer 
 points\footnote{i.e., construction of analytic continuation
  involves the Carlson theorem, see, e.g.,
 \cite{AAR}, Theorem 2.8.1}" of the Meixner--Pollaszek
 polynomials.

%%%%%%%%%%%%%%%%%%%%%%%%%%%%%%%%%%%%%%%%%%%%%%%%%%%%%

%%%%%%%%%%%%%%%%%%%%%%%%%%%%%%%%%%%%%%%%%%%%%%%%%%%%%%%%%%%%%%%%%%%%%%%
%%%%%%%%%%%%%%%%%%%%%%%%%%%%%%%%%%%%%%%%%%%%%%%%%%%%%%%%%%%%%%%%%%%%%%%%%%
%%%%%%%%%%%%%%%%%%%%%%%%%%%%%%%%%%%%%%%%%%%%%%%%%%%%%%%%%%%%%%%%%%%%%%%%%%
%%%%%%%%%%%%%%%%%%%%%%%%%%%%%%%%%%%%%%%%%%%%%%%%%%%%%%%%%%%%%%%%%%%%%%%%%%%%

\section{Preliminaries}

\COUNTERS

{\bf\punct Imaginary shift in $L^2$.}
We say that a function is holomorphic in a closed strip 
$|\Im s|\le \alpha$ if it is holomorphic in a larger strip 
$|\Im s|<\alpha+\delta.$

\begin{lemma}
\label{l:l-2.1}
Let $H\subset L^2(\R)$ be the subspace in $L^2(\R)$ consisting of functions $f(s)$ admitting
holomorphic continuation to the strip
$|\Im s|\le 1$ and satisfying the condition 
$|f(s)|=O(s^{-1/2-\epsilon})$ in this strip.
The operators 
$$
T_+f(s)=f(s+i), \qquad T_-f(s)=f(s-i)
$$
defined on $H$ 
are symmetric in $L^2(\R)$.
\end{lemma}

{\sc Proof.}
$$
\int_{-\infty}^\infty f(s+i)\,\ov{g(\ov s)}\,ds=
\int_{i-\infty}^{i+\infty} f(t) \ov{g(\ov t+i)}\,dt=
\int_{-\infty}^{\infty} f(t) \ov{g(\ov t+i)}\,dt
.$$

{\bf\punct Lemma on symmetry.}
Now let $\mu(s)$, $\nu(s)$ be the same as above, see
 (\ref{eq:mu})--(\ref{eq:nu}). Therefore the weight
$w(s)$ is 
$$
w(s)=\mu(s)\nu(s)=\frac 1{2\pi}
e^{2cs}
\frac{\prod_{k=1}^m \Gamma(a_k+ is)\Gamma(\ov a_k+- is)}
{\prod_{l=1}^n\Gamma(b_l+is)\Gamma(\ov b_l-is)}
.
$$
For real $s$ we can represent $w(s)$ in the form
$$
w(s)=\frac 1{2\pi}
e^{2cs}
\left|\frac{\prod_{k=1}^m \Gamma(a_k+ is)}{\prod_{l=1}^n\Gamma(b_l+is)}
\right|^2
.
$$

Let $A(s)$, $B(s)$ be as above
\begin{align*}
A(s)&:=\frac{\nu(s+i)}{\nu(s)}= e^{-i c}\frac{\prod_{k=1}^m 
(\ov a_k- is)}{\prod_{l=1}^n(\ov b_l-is)}
\\
B(s)&:=\frac{\mu(s-i)}{\mu(s)}= e^{ic}\frac{\prod_{k=1}^m 
(a_k+ is)}{\prod_{l=1}^n(b_l+is)}
.
\end{align*}

By (\ref{eq:gamma-is}), we have the following asymptotics of $w(s)$
 in any strip $|\Im s|<\alpha $
\begin{equation}
w(s)\sim \Psi(s):=\mathrm{const}\cdot |s|^{\sum(2 \Re a_k-1)-\sum(2 \Re b_l-1)}
 \exp\bigl(2cs+(n-m)\pi s\bigr),
\end{equation}
as 
$s\to\infty$.
We say that a function $f$ is {\it $w$-decreasing} in
 a strip $|\Im s|\le\alpha$
 if 
 $$
 f(s)=O\left(\Psi(s)^{-1/2}
s^{-m-1/2-\epsilon}\right), \qquad s\to\infty 
. $$
 
 This condition provides 
 $$
 f(s+i\beta),\,\, A(s)f(s+i\beta),\,\, B(s)f(s+i\beta)\in L^2(\R, w(s)\,ds)
 $$
  for all
 $\beta$ satisfying $|\beta|\le \alpha$.
 Denote by $\cH[w]$ the space of all functions holomorphic in the strip
 $|\Im s|\le 1$ and $w$-decreasing in this strip. 
 \hfill $\square$

\begin{lemma}
\label{l:l-2.2}
Let
$$
\Re a_j>0
$$
for all $j$.
The operator
$$
Rf(s)=A(s) f(s+i)
$$
defined on the domain $\cH[w]$
is symmetric in $L^2(\R,w(s)\,ds)$.
\end{lemma}

{\sc Proof.}
We verify the identity
$\la Rf,g\ra=\la f, Rg\ra$ for $f$, $g\in \cH[w]$:
\begin{multline*}
\int_{-\infty}^\infty
\frac{\nu(s+i)}{\nu(s)}f(s+i)\,\ov{g(\ov s)} \mu(s)\nu(s)\,ds
=
\int_{-\infty}^\infty
f(s+i)\,\ov{g(\ov s)} \mu(s)\nu(s+i)\,ds
=\\=
\int_{i-\infty}^{i+\infty}
f(s)\,\ov{g(\ov s+i)} \mu(s-i)\nu(s)\,ds=
\int_{-\infty}^\infty
f(s)\,\ov{g(\ov s+i)} \mu(s-i)\nu(s)\,ds=
\\=
\int_{-\infty}^\infty
f(s)\,\frac{\mu(s-i)}{\mu(s)} \ov{g(\ov s+i)}  \mu(s)\nu(s)\,ds
=
\int_{-\infty}^\infty
f(s)\,\ov{\frac{\nu(\ov s+i)}{\nu(\ov s)} g(\ov s+i)}  \mu(s)\nu(s)\,ds
.
\end{multline*}
The condition $\Re a_j>0$ provides absence of poles of
$\nu(s+i)\mu(s)$ in the strip $0<\Im s<1$.

\begin{corollary}
Under the same conditions the operator
$$
\cL f(s)=A(s)f(s+i)-(A(s)+B(s)) f(s)+ B(s) f(s-i)
$$
is symmetric on the subspace $\cH[w]\subset L^2(\R,w(s)\,ds)$.
\end{corollary}

\sm

%%%%%%%%%%%%%%%%%%%%%%%%%%%%%%%%%%%%%%%%%%%%%%%%%%%%%%%%%%%%%%%%%%%%%%%%%%%

{\bf\punct Change of a weight.}
Let $w_2(s)=\tau(s)\ov{\tau(\ov s))} w_1(s)$. 
Then the operator
$$
Hf(s)=\tau(s) f(s)
$$
is a unitary operator $L^2(\R,w_2(s))\to L^2(\R,w_1(s))$.
Evidently,
$$
H^{-1} T_+ H f(s)=
\frac {\tau(s+i)}{\tau(s)} T_+,\qquad
H^{-1} T_- H f(s)=
\frac {\tau(s-i)}{\tau(s)} T_+
.$$

%%%%%%%%%%%%%%%%%%%%%%%%%%%%%%%%%%%%%%%%%%%%%%%%%%%%%%%%%%

\sm

{\bf\punct Operators in $L^2(\R)$.}

\begin{lemma}
Let an operator
$$
R f(s)=L(s)f(s+i)
$$
be formally symmetric in $L^2(\R, ds)$. Then
\begin{equation}
L(s)=\ov{L(\ov s-i)}.
\label{eq:funct-eq}
\end{equation}
\end{lemma}

This is straightforward. 

Note that if $L(s)$ satisfy (\ref{eq:funct-eq}), then
$L(s)^{-1}$  satisfy the same condition. Also, if $L_1(s)$, $L_2(s)$
satisfy (\ref{eq:funct-eq}), then $L_1(s)L_2(s)$ satisfy (\ref{eq:funct-eq}).

Obvious solutions are
\begin{align*}
L(s)&=i/2+s,
\\
L(s)&=(i/2+ia+s)(i/2-ia+s),
\\
L(s)&=h(e^{2\pi s})
.
\end{align*}

%%%%%%%%%%%%%%%%%%%%%%%%%%%%%%%%%%%%%%%%%%%%%%%%%%%%%%%%%%%%%%%%%%%%%
%%%%%%%%%%%%%%%%%%%%%%%%%%%%%%%%%%%%%%%%%%%%%%%%%%%%%%%%%%%%%%%%%%%%%%%%%
%%%%%%%%%%%%%%%%%%%%%%%%%%%%%%%%%%%%%%%%%%%%%%%%%%%%%%%%%%%%%%%%%%%%%%%%%%%%%%

\section{The Kontorovich--Lebedev transform}

\COUNTERS

{\bf\punct Difference operator.}
Now $\mu(s)=\Gamma(is)$, $w(s)=|\Gamma(is)|^{-2}$.
 We consider the space of {\it even} functions,
$f(s)=f(-s)$,
the inner product is given by
$$
\la f,g\ra=\frac 2 {\pi}\int_{-\infty}^\infty
f(s)\ov{g(\ov s)} \frac {ds}{|\Gamma(is)|^2}
=
\frac2{\pi^2}
\int_{-\infty}^\infty
f(s)\ov{g(\ov s)} s\sinh (\pi s)\,ds
.
$$
We consider a difference operator
$\cL$ given by
\begin{equation}
\cL f(s)= \frac 1{is}\bigl(f(s+i)-f(s-i)\bigr)
\label{eq:cL-KL}
\end{equation}
defined on the subspace $\cH[w]\subset L^2(\R_+,|\Gamma(is)|^{-2}ds)$ 

\begin{lemma}
\label{l:sa-KL}
 The operator $\cL$ is essentially self-adjoint.
\end{lemma}

The spectral decomposition is given by the inverse Kontorovich--Lebedev transform, see the next 
subsection.

\smallskip

%%%%%%%%%%%%%%%%%%%%%%%%%%%%%%%%%%%%%%%%%%%%%%%%%%%%%%%%%%%%%%%%%%%%%%%%%%55

{\bf\punct The Kontorovich--Lebedev transform. Preliminaries.}
The {\it Macdonald functions} $K_\nu(z)$ are solutions of the  modified
Bessel differential equation  (see \cite{HTF2},7.2(11)), i.e. the equation
$$
\left(z\frac d{dz}\right)^2 g(z)-z^2 g(z)=-\nu^2 g(z)
.
$$
They are defined by (see \cite{HTF2}, 7.2(13)),
$$
K_\nu(z)=
\frac\pi{\sin(\nu \pi)}
(I_{-\nu}(z)-I_{\nu}(z))
,$$
where $I_\nu(z)$ are the modified Bessel functions,
$$
I_\nu(z)=e^{-i\nu\pi/2}
J_\nu(ze^{i\pi/2})=
\sum_{m=0}^\infty \frac{(z/2)^{2m+\nu}}{m!\,\Gamma(m+\nu+1)}
.
$$
For each $z\ne 0$ the function $h_z(\nu):=K_\nu(z)$ is an entire function of the variable
$\nu$, 
$$
K_\nu(z)=K_{-\nu}(z)
.
$$
For positive $z\in\R$ and $\nu\in i\R$  values of $K_\nu(z)$ are real.

Below we use two identities (see \cite{Wat}, (3.71.1)--(3.71.2))
\begin{align}
K_{\nu-1}(z)-K_{\nu+1}(z)=-\frac{2\nu}z K_\nu(z)
\label{eq:W1}
,
\\
K_{\nu-1}(z)+K_{\nu+1}(z)=-\frac d{dz}K_\nu(z)
\label{eq:W2}
.
\end{align}

The Kontorovich--Lebedev transform \cite{KoLe}, \cite{Leb}, Section 6.5, 
\cite{YaLu} is given by%
\footnote{Here and below we understand  integral operators
in the sense of the kernel theorem, see, e.g., \cite{Her}, Section 5.2. 
However, for the Kontorovich--Lebedev transform and the Wimp transform discussed
below conditions of literal validness of formulas are well investigated.}
the formula 
\begin{equation}
\frK g(s)=\int_0^\infty K_{is} (x) g(x)\frac{dx}x
.
\end{equation}

The inverse transform is 
\begin{equation}
\frK^{-1} f(x)=\frac 2\pi \int_0^\infty 
f(s) K_{is}(x) \frac{ds}{|\Gamma(is)|^2}
.
\end{equation}
The Kontorovich--Lebedev transform is a unitary operator
$$
L^2(\R_+,x^{-1}dx)\to L^2(\R_+, 2\pi^{-1}|\Gamma(is)|^{-2}ds)
.
$$

%%%%%%%%%%%%%%%%%%%%%%%%%%%%%%%%%%%%%%%%%%%%%%%%%

{\bf\punct The statement.}

\begin{theorem}
The Kontorovich-Lebedev transform provides a unitary equivalence
between the operator
$$
Pg(x)=\frac 2x g(x)
$$
in $L^2(\R_+,x^{-1}dx)$ and the operator
$\cL$ given by {\rm(\ref{eq:cL-KL})}.
\end{theorem}

{\sc Proof.} We use (\ref{eq:W1}),
\begin{multline*}
\int_0^\infty \frac 2x g(x)\cdot K_{is}(x)\frac {dx}x=
\int_0^\infty  g(x)\cdot \frac{2 K_{is}(x)}{x}\frac {dx}x
=\\=
\int_0^\infty g(x) \frac 1{is} \bigl(K_{i(s+i)}(x)-K_{i(s-i)}(x)\bigr)\frac {dx} x
=\frac 1{is}\bigl(\frK g(s+i)-\frK g(s-i)\bigr)
.
\end{multline*}
This proves the statement. 

\sm

{\sc Remark.} However,  Lemma \ref{l:sa-KL} in this moment is not proved, 
it a special case of Lemma \ref{l:sa-W} proved below.

%%%%%%%%%%%%%%%%%%%%%%%%%%%%%%%%%%%%%%%%%%%%%%%%%%%

\sm

{\bf\punct An additional remark.} Applying
(\ref{eq:W2}), we get the following statement

\begin{proposition}
The Kontorovich--Lebedev transform send
the operator
$$
Qg(x)=\left(\frac {d}{dx}-\frac 1x\right) g(x)
$$
to the operator
$$
\cM f(s)=\frac 12\bigl(f(s+i)-f(s-i)\bigr)
.
$$
\end{proposition}
 
Therefore, we can evaluate the image of any operator
$x^{-m}\frac {d^n}{dx^n}$ under the Kontorovich--Lebedev transform.

\section{The Wimp transform}

\COUNTERS

{\bf\punct Difference problem.}
Now 
$
\mu(s)=\frac{\Gamma(1/2-\rho+is)}{\Gamma(2is)}
$,
we consider the space of {\it even} functions on $\R$ with inner product
$$
\la f,g\ra =
\frac 1{4\pi}
\int_{-\infty}^{\infty} f(s)\ov{g(\ov{s})} 
\left|\frac{\Gamma(1/2-\rho+is)}{\Gamma(2is)}\right|^2
ds
.
$$
We consider the following difference operator 
\begin{multline}
\cL f(s)=\frac{1-\rho-is}{(-2is)(1-2is)}f(s+i)-
\\-
\left( \frac{1-\rho-is}{(-2is)(1-2is)}+ \frac{1-\rho+is}{(+2is)(1+2is)}\right) f(s)
+ \frac{1-\rho+is}{(+2is)(1+2is)}f(s-i)
.
\label{eq:cL-wimp}
\end{multline}
As above, this operator is defined on the subspace $\cH[w]\subset L^2(\R,w(s)ds)$.

\begin{lemma}
\label{l:sa-W}
Let $\rho<1/2$. Then 
the operator $\cL$ is essentially self-adjoint.
\end{lemma}

%%%%%%%%%%%%%%%%%%%%%%%%%%%%%%%%%%%%%%%%%%%%%%%%%%%%%%%%%%%%%%%%%%%%%%%%

{\bf\punct Whittaker functions and the Wimp transform. Preliminaries.}
{\it Whittaker functions} $W_{\rho,\sigma}(z)$ are versions of confluent hypergeometric
functions. They are solutions of the Whittaker equation (see \cite{HTF1}, 6.1 (4))
\begin{equation}
\left(x^2\frac {d^2}{dx^2} -\frac{x^2}4+\rho x \right) f(x)
=(\sigma^2-1/4) f(x)
.
\label{eq:whitt}
\end{equation}
The explicit expression is
\begin{multline}
W_{\rho,\sigma} (x)=
e^{-x/2} 
\Bigl(
\frac{\Gamma(-2\sigma) x^{1/2+\sigma}}{\Gamma(1/2-\rho-\sigma)}\,{}_1F_1\left[
\begin{matrix}1/2-\rho+\sigma\\ 1+ 2\sigma \end{matrix}; x\right]
+\\+
\frac{\Gamma(2\sigma) x^{1/2-\sigma}}{\Gamma(1/2-\rho+\sigma)}\,{}_1F_1\left[
\begin{matrix}1/2-\rho-\sigma\\ 1- 2\sigma \end{matrix}; x\right]
\Bigr)
\label{eq:W-F11}
,\end{multline}
see %\cite{HTF1}, 6.9(4), 6.5(7).
\cite{Sla}, (1.9.10).
There are the following integral representations 
(see \cite{HTF1}, 6.11(18), \cite{PBM1}, 2.3.6.9),
\begin{equation}
W_{\rho,\sigma}(x)=
\frac{e^{-x/2}x^{\rho}}{\Gamma(1/2-\rho+\sigma)}
\int_0^\infty 
e^{-xt} t^{-1/2-\rho+\sigma}(1+t)^{-1/2+\rho+\sigma}dt
.
\label{eq:kummer}
\end{equation}
 and the Barnes representation (see \cite{PBM3}, 8.4.44.3, \cite{Sla},(3.5.16)),
\begin{multline}
W_{\rho,\sigma} (z)=
\frac 
{ e^{-x/2}}
{2\pi\Gamma(1/2-\rho-\sigma) \Gamma(1/2-\rho+\sigma)}
\times\\ \times
\int_{-\infty}^\infty
\Gamma(it+1/2+\sigma) \Gamma(it+1/2-\sigma)\Gamma(-\rho-it)
x^{-it}dt
.
\label{eq:W-barnes}
\end{multline}

{\sc Remark.} If $\rho\in \R$,  $\sigma\in i\R$, $x>0$, then
$W_{\rho,\sigma}(x)$ is real. This follows from (\ref{eq:W-F11}).
\hfill $\square$

\sm

Fix real $\rho<1/2$. The {\it Wimp transform} $\frW_\rho$ is 
the integral operator given by
$$
\frW_\rho g(s)=\int_0^\infty g(x) W_{\rho,is}(x) \frac{dx}{x^2}
$$
(see \cite{Wim}, \cite{YaLu}).
The inverse transform is
$$
\frW_\rho^{-1} f(x)=\frac 1{2\pi} \int_0^\infty f(x) W_{\rho,is}(x) 
\left|\frac{\Gamma(1/2-\rho+is)}{\Gamma(2is)}\right|^2
ds
.
$$
The Wimp transform is a unitary operator
$$
L^2(\R_+,x^{-2}dx)\to L^2\left(\R_+, 
\frac1{2\pi}\left|\frac{\Gamma(1/2-\rho+is)}{\Gamma(2is)}\right|^2
ds  \right)
.
$$

{\sc Remark.} This theorem can be obtained by writing of explicit spectral decomposition
of the differential operator (\ref{eq:whitt}) as it 
is explained in \cite{Dunford}, Chapter XIII.
\hfill $\square$

\sm

The Macdonald function $K_\nu$ admits the following expression in the terms
of Whittaker functions:
$$
K_\nu(x)=\sqrt{\frac\pi{2x}} W_{0,\nu}(x)
.
$$
Therefore the Kontorovich--Lebedev transform is a special case of
Wimp transforms.

%%%%%%%%%%%%%%%%%%%%%%%%%%%%%%%%%%%%%%%%%%%%%%%%%%%%%%%%%%%%%%%%%%%

\sm

{\bf\punct The statement.}

\begin{theorem}
The Wimp transform send the operator
\begin{equation}
Rg(x)=x^{-1} g(x)
\label{eq:}
\end{equation}
to the difference  operator $\cL$ defined by
{\rm(\ref{eq:cL-wimp})}.
\end{theorem}

Theorem is a corollary of the following lemma.

\begin{lemma}
The Whittaker functions satisfy the difference equation
\begin{multline}
\frac{1-\rho-\sigma}{(-2\sigma)(1-2\sigma) }( W_{\rho,\sigma-1}(x)
-W_{\rho,\sigma})+
\frac{1-\rho+\sigma}{(2\sigma)(1+2\sigma) }( W_{\rho,\sigma+1}(x)
-W_{\rho,\sigma})=
\\=
\frac 1x W_{\rho,\sigma}(x)
\label{eq:diff-eq-W}
.
\end{multline}
\end{lemma}
 
{\sc Proof.} We use the Barnes integral (\ref{eq:W-barnes}).
We multiply both sides of (\ref{eq:diff-eq-W}) by $e^{x/2}$
and pass to Mellin transforms 
(see  below (\ref{eq:mellin})--(\ref{eq:mellin-inverse})). 
Denote by $h(t)$ the Mellin transform of $e^{x/2}W_{\rho,\sigma}(x)$, i.e.,
$$
h(t)=
\frac
{\Gamma(it+1/2+\sigma) \Gamma(it+1/2-\sigma)\Gamma(-\rho-it)}
{\Gamma(1/2-\rho-\sigma) \Gamma(1/2-\rho+\sigma)}
.
$$
The Mellin transforms of $e^{x/2}W_{\rho,\sigma\pm1 }(x)$ are
$\gamma_{\pm}(t)h(t)$, where 
$$
\gamma_\pm(t)=\frac{(-1/2-\rho\mp\sigma)(t+1/2\pm\sigma)}
{(1/2-\rho\pm\sigma)(t-1/2\mp\sigma)}
.
$$
In the left-hand side we get 
 \begin{multline*}
h(s)\cdot\left\{ \frac{1-\rho-\sigma}{(-2\sigma)(1-2\sigma) }(\gamma_-(t)-1)
+
\frac{1-\rho+\sigma}{(2\sigma)(1+2\sigma) }(\gamma_+(t)-1)\right\}
=\\
=h(s)\frac{-t-\rho}{(t-1/2-\sigma)(t-1/2+\sigma)}
=\\=
\frac
{\Gamma(it-1/2+\sigma) \Gamma(it-1/2-\sigma)\Gamma(-\rho-it+1)}
{\Gamma(1/2-\rho-\sigma) \Gamma(1/2-\rho+\sigma)}=
h(s+i)
.
 \end{multline*}
 Shift of a Mellin transform by $i$ is equivalent to multiplication of the original by
 $1/x$.
 \hfill $\square$
 
 %%%%%%%%%%%%%%%%%%%%%%%%%%%%%%%%%%%%%%%%%%%%
 
 \sm
 
{\bf\punct Proof of self-adjointness.} 
The space $C^{\infty}_c(\R_+)$ of smooth functions 
with compact support on
 $(0,\infty)$ is a domain of essential self-adjointness
of the operator (\ref{eq:}). It is sufficient to prove 
the following lemma.

\begin{lemma}
\label{l:W-image}
$\frW_\rho\bigl(C^{\infty}_c(\R_+)\bigr)\subset\cH[w]$.
\end{lemma}

\begin{lemma}
\label{l:compact}
Fix $\rho<1/2$.
For $(\sigma,x)$ ranging in a domain 
\begin{equation}
 |\Re\sigma|\le 1 , \qquad 0<c\le x\le C< \infty
\label{eq:compact}
\end{equation}
the following uniform estimate holds
\begin{equation}
|W_{\rho,\sigma}(x)|=O\left(  e^{\pi |\Im\sigma|/2} 
|\Im\sigma|^{\rho+1}\right)
\label{eq:W-estimate}
.
\end{equation}
\end{lemma}

{\sc Proof of Lemma \ref{l:compact}.}
The integral  formula (\ref{eq:kummer})  converges
if $\Re\sigma>\rho-1/2$ and admits the holomorphic
to the whole plane $\sigma\in\C$.

The statement is very simple if $\rho<-1/2$ (the integral in
(\ref{eq:kummer}) is bounded and the desired estimate is obtained from
an estimate of a pre-integral factor. But we wish to cover also
the interval $-1/2<\rho<1/2$.

Fix $A>B>1$.
Represent $1$ as $1=\phi(t)+\psi(t)$, where $\phi$, $\psi(t)\ge 0$
are smooth nonnegative on $\R_+$, $\psi(t)=0$ for $t<A$, and
$\phi=0$ for  $t>B$.
We write the integral in (\ref{eq:kummer}) as
\begin{multline}
\int_0^\infty=
\int_0^A
e^{-xt} t^{-1/2-\rho+\sigma}(1+t)^{-1/2+\rho+\sigma}\phi(t)\,dt
+\\+
\int_B^\infty 
e^{-xt} t^{-1/2-\rho+\sigma}(1+t)^{-1/2+\rho+\sigma}\psi(t)\,dt
.
\label{eq:two-summands}
\end{multline}
The second summand is uniformly bounded in our domain
(\ref{eq:compact}),  the integrand is dominated  by
$$
e^{-ct} t^{1/2-\rho}(1+t)^{1/2+\rho}
$$

Next, we represent the first summand of (\ref{eq:two-summands}) as
$$
\int_0^A=
\int_0^A
 t^{-1/2-\rho+\sigma} 
\left(e^{-xt} (1+t)^{-1/2+\rho+\sigma}\phi(t)-1\right)\,dt
+
\int_0^A  t^{-1/2-\rho+\sigma} \,dt
$$
Denote by $Q(t,x,\sigma)$ the first integrand.
Then $|Q(t,x,\sigma)|$ depend on $t$, $x$, $\Re\sigma$,
 these variables
range in a compact set, the function $Q$ is continuous on this set.
Therefore
first summand is uniformly bounded in (\ref{eq:compact}),
 the second 
summand is uniformly bounded in (\ref{eq:compact}) 
outside a neighborhood of 
$\sigma=\rho-1/2$.

Thus $\int_0^\infty$ is uniformly bounded in 
in (\ref{eq:compact}) 
outside a neighborhood of 
$\sigma=\rho-1/2$.

Next, we  multiply the integral  (\ref{eq:kummer}) 
 by the pre-integral factor 
$\frac{e^{-x/2}x^{\rho}}{\Gamma(1/2-\rho+\sigma)}$.
Since $\Re(1/2-\rho+\sigma)\in (-1/2-\rho,3/2-\rho)$, we have
$$
\Gamma(1/2-\rho+\sigma)^{-1}
= O(e^{\pi |\Im\sigma|/2} |\Im \sigma|^{\rho+1}), \qquad 
|\Im \sigma|\to\infty
$$
and we get (\ref{eq:W-estimate}). 
\hfill$\square$

\sm

%%%%%%%%%%%%%%%%%%%%%%%%%%%%%%%%%%%%%%%%%%%%%%%%%%%%%%%%%%%%%

{\sc Proof of Lemma \ref{l:W-image}.}
By Lemma \ref{l:compact}, 
 for a function $f\in C^\infty_c(\R_+)$ with compact support,
we have 
\begin{equation}
|\frW_\rho f(s)| \le C\cdot e^{\pi |\Re s|/2} |\Re s|^{\rho+1/2}
.
\label{eq:estimate}
\end{equation}
Next, we use (\ref{eq:whitt}),
\begin{multline*}
-(1/4+s^2) \frW_\rho f(s)=\int_0^{\infty} (-1/4-s^2) W_{\rho,is}(x))\cdot f(x)\frac{dx}x=
\\=
\int_0^{\infty}\left(x^2\frac {d^2}{dx^2} -\frac{x^2}4+\rho x \right)
W_{\rho,is}(x)\cdot f(x)\frac{dx}x
=\\=\int_0^\infty W_{\rho,is}(x) \cdot
\left[ x\frac{d^2}{dx^2}(x f(x))-\frac14 x^2 f(x)+\rho x f(x)\right]
\frac{dx}x
.
\end{multline*}
We apply (\ref{eq:estimate}) for the function
 in square brackets and get
$$
|\frW_\rho f(s)| \le C\cdot
   (s^2+1/4)^{-1}\cdot e^{\pi |\Re s|/2} |\Re s|^{\rho+1}
$$
%Repeating  the same trick, we come to
and $\frW_\rho f(s)\in \cH[w]$.
\hfill $\square$

%%%%%%%%%%%%%%%%%%%%%%%%%%%%%%%%%%%%%%%%%%%%%%%%%%%%%%%%%%%
%%%%%%%%%%%%%%%%%%%%%%%%%%%%%%%%%%%%%%%%%%%%%%%%%%%%%%%%%%
%%%%%%%%%%%%%%%%%%%%%%%%%%%%%%%%%%%%%%%%%%%%%%%%%%%%%%%%%%%%
%%%%%%%%%%%%%%%%%%%%%%%%%%%%%%%%%%%%%%%%%%%%%%%%%%%%%%%%%%%%

\section{The Vilenkin transform}

\COUNTERS

{\bf\punct Difference problem.} Fix $\alpha>0$, $\phi>0$.
We consider the weight
$$
w(t)=\frac 1{2\pi}|\Gamma(\alpha/2+it)|^2 e^{\pi t}
,$$
 the corresponding space $L^2(\R,w(t)\,dt)$, and the difference operator
$$
\cL f(t)=i(\alpha/2+it) f(t-i)+ 2 \cosh \phi \,t h(t)- i(\alpha/2-it) f(t+i)
.
$$
This operator differs from (\ref{eq:mu})-(\ref{eq:difference-operator}),
 but 
it is symmetric (proof is the same as in  Lemma \ref{l:l-2.2}).

\begin{theorem}
The operator $\cL$ is essentially self-adjoint on the space $\cH[w]$.
\label{th:vilenkin-1}
\end{theorem}

%%%%%%%%%%%%%%%%%%%%%%%%%%%%%%%%%%%%%%%%%%%%%%%%%%%%%%%

\sm

{\bf \punct The Vilenkin transform.}

\begin{theorem}
\label{th:vilenkin-2}
 The Vilenkin transform
$$
\frV_\alpha g(t)
= (1-e^{-2\phi})^{\alpha/2} e^{-\phi it}
\int\limits_{-\infty}^\infty
g(s)\,\,
{}_2F_1 \left[
\begin{matrix} \alpha/2-is,\alpha/2+it\\ \alpha\end{matrix};1-e^{-2\phi}  \right]
w(s)\,ds
.
$$
is a unitary operator $L^2(\R,w(s)\,ds)\to L^2(\R,w(s)\,ds)$.
\end{theorem}

This  is  a minor  modification 
of Vilenkin \cite{Vil}, \S 7.4, see also \cite{VK}, 7.7.7.

Since the operator $\frV$ is unitary, the inversion formula is
$$
\frV_\alpha^{-1} f(s)
= (1-e^{-2\phi})^{\alpha/2} 
\int\limits_{-\infty}^\infty
f(t)\,\,
{}_2F_1 \left[
\begin{matrix} \alpha/2+is,\alpha/2-it\\ \alpha\end{matrix};1-e^{-2\phi}  \right]
e^{\phi it}
w(t)\,dt
.
$$

\begin{theorem}
\label{th:vilenkin-3}
The inverse Vilenkin transform $\frV_\alpha^{-1}$ send the operator $\cL$
to the operator 
$$
\cN f(s)=2s\sinh \phi f(s)
.
$$
\end{theorem}

To prove these statements,
 we decompose the Vilenkin transform as a product
of 3 simple transformations, see below formula (\ref{eq:JTJ}).

\sm

%%%%%%%%%%%%%%%%%%%%%%%%%%%%%%%%%%%%%%%%%%%%%%%%%%%%%%%%%%%%%

{\bf\punct Highest weight representations of $\SL_2(\R)$.}
The group $\SL(2,\R)$ is the group of $2\times 2$ real matrices
$g=\begin{pmatrix}a&b\\c&d  \end{pmatrix}$ with $\det g=1$.
Denote by $\Pi$ the half-plane $\Im z>0$.

Fix $\alpha>0$. Consider the Hilbert space $H_\alpha$ of holomorphic functions
on $\Pi$ determined by reproducing kernel
(see, e.g., \cite{Ner-gauss}, Section 7.1)
$$
K(z,u)=\left(\frac{z-\ov u}{2i}\right)^{-\alpha}
.$$
In other words, denote $\Psi_a(z):=K(z,a)$. Then for any $F\in H_\alpha$
we have
\begin{equation}
\la F,  \Psi_a\ra =F(a)
\label{eq:reproducing}
.
\end{equation}
For $\alpha>1$ the inner product in $H_\alpha$ admits
the following
integral representation
$$
\la F, G\ra=\const(\alpha) 
\int_{\Pi} F(z)\ov{G(z)} (\Im z)^{\alpha-2} dz\,d\ov z
.
$$
Consider the following operators
in $H_\alpha$
$$
T_\alpha\begin{pmatrix}a&b\\c&d\end{pmatrix}
F(z)=F\left(\frac{b+zd}{a+zc} \right)(a+zc)^{-\alpha}
.
$$
The function $(a+zc)^{-\alpha}$ is multi-valued. We choose arbitrary branch of
this function on $\Pi$.
Then operators $T_\alpha(g)$ are unitary and satisfy the condition
$$
T_\alpha(g_1)T_\alpha(g_2)=\lambda(g_1,g_2)
T_\alpha(g_1 g_2)
,$$ 
where $\lambda(g_1,g_2)\in \C$. Thus we get a projective unitary representation
of $\SL_2(\R)$, such representations are called {\it highest weight representations}.

%%%%%%%%%%%%%%%%%%%%%%%%%%%%%%%%%%%%%%%%%%%%%%%%%%%%%%%%%%%%%%%

\sm

{\bf\punct The Mellin transform. Preliminaries.} See, e.g., \cite{Tit}.
 For a function $f$ on $\R_+$ we define a Mellin transform
$\frM f(s)$ as
\begin{equation}
\frM f(s)=\int_0^\infty f(x) x^{is-1} dx
.
\label{eq:mellin}
\end{equation}
The inverse transform is given by 
\begin{equation}
\frM^{-1} g(x)=\frac 1{2\pi} \int_{-\infty}^\infty g(s) x^{-is} ds
.
\label{eq:mellin-inverse}
\end{equation}
The Mellin transform is a unitary operator
$L^2(\R_+,x^{-1}dx)\to L^2(\R,\frac 1{2\pi} ds)$.

Notice, that changing variable $x=e^t$ in (\ref{eq:mellin}), we come to the usual Fourier
transform.

%%%%%%%%%%%%%%%%%%%%%%%%%%%%%%%%%%%%%%%%%%%%%%%%%%%%%%%%%%%%%

\sm

{\bf\punct Spectral decomposition of dilatation operators.}
Consider a one-parametric subgroup $A\simeq \R_+^\times$
 in $\SL_2(\R)$ consisting of matrices
of the form
$D(a)=\begin{pmatrix} a^{1/2}&0\\0&a^{-1/2}
\end{pmatrix}$, where $a>0$.
The subgroup $A$ acts in the space $H_\alpha$ by the transformations
\begin{equation}
T_\alpha(D(a))=f(z)
=f(a^{-1}z)a^{-\alpha/2}
\label{eq:D1}
.
\end{equation}
Next, consider the measure $d\mu(s)$ on $\R$ given by
$$
\mu(s)\,ds:=\frac 1{2\pi\Gamma(\alpha)} |\Gamma(\alpha/2+is)|^2\,ds
$$ 
and the action of the same group in the space
$L^2(\R,\mu(s)\,ds)$ given by the formula
\begin{equation}
\tau_\alpha(D(a))f(s)=f(s)a^{is}
\label{eq:D2}
.
\end{equation}
Consider the operator
 $J:L^2(\R,\mu(s)\,ds)\to H_\alpha$
given by
\begin{equation}
F(z)=J_\alpha f(z)=\frac{2^{\alpha}}{2\pi\Gamma(\alpha)}
\int_{-\infty}^\infty f(s) \left(\frac zi\right)^{-\alpha/2-is} |\Gamma(\alpha/2+is)|^2ds
\label{eq:mellin-twist}
,
\end{equation}
we choose a branch of $(z/i)^{-\alpha/2-is}=e^{-(\alpha/2+is)\ln(z/i)}$
such  that $\ln z/i$ is real for $z=ip$, $p>0$.

Therefore $F(ip)(p)^{\alpha/2}$ is the inverse Mellin transform of 
$\frac {2^\alpha}{\Gamma(\alpha)}f(s)|\Gamma(\alpha/2+is)|^2$.

Applying the direct Mellin transform, we get
\begin{equation}
f(s)\cdot \frac{2^{\alpha}}{\Gamma(\alpha)}|\Gamma(\alpha/2+is)|^2 
=\int_0^\infty F(ip)\,p^{\alpha/2+is-1} dp
.
\end{equation}

\begin{proposition}
The transform $J_\alpha$ is a unitary operator
$$J_\alpha:L^2(\R,\mu(s)\,ds)\to H_\alpha$$
intertwining actions {\rm (\ref{eq:D1})} and {\rm (\ref{eq:D2})}.
\end{proposition}

{\sc Proof.} A verification of
$$
J\circ \tau_\alpha(D(a))=T_\alpha (D(a))\circ J
$$
is straightforward.
 Next,
$$
T_\alpha(D(a)) \Psi_i=a^{\alpha/2} \Psi_{a i} 
.
$$
By (\ref{eq:reproducing}), the system of vectors $\Psi_{a i}$, where $a>0$
 is total
in the Hilbert space $H_\alpha$.

Next, we consider functions 
$$\Phi_a:=a^{\alpha/2+is}$$
in $L^2(\R,\mu(s)\,ds$. Then 
$$
\tau_\alpha(D(a)) \Phi_1=\Phi_a\cdot a^{\alpha/2}
.
$$ 
To prove unitarity, it is sufficient to show (see, e.g., \cite{Ner-gauss}, Theorem 7.1.4) that
\begin{align}
 J\Phi_a&=\Psi_{ia}
, \label{eq:JPhi}
 \\
\la\Phi_a,\Phi_b \ra_{L^2}&=
\la \Psi_{ia},\Psi_{ib}\ra_{V_{\alpha}}=\left(\frac{a+b}2\right)^{-\alpha}
.
\label{eq:PhiPhi}
 \end{align}
 First, note that 
 $$
 \int_0^\infty (1+x)^{-\alpha} x^{-is}dx=\mathrm B(is,\alpha-is)=
 \frac{\Gamma(is)\Gamma(\alpha-is)}{\Gamma(\alpha)}
. $$
 Applying the inversion formula for the Mellin transform, we
 get (\cite{PBM3}, 8.5.2.5)
 $$
 \frac 1{2\pi}\int_{-\infty}^\infty \Gamma(is)\Gamma(\alpha-is) x^{is-1}ds=
   {\Gamma(\alpha)} (1+x)^{-\alpha}
 .
 $$
 Both formulas (\ref{eq:JPhi})--(\ref{eq:PhiPhi}) are reduced to the latter integral.
 \hfill $\square$
 
 %%%%%%%%%%%%%%%%%%%%%%%%%%%%%%%%%%%%%%%%%%%%%%%%%%%%%%%%%%%%
 
\sm 
 
 {\bf\punct Calculation. Proof of Theorem \ref{th:vilenkin-2}}
 Set
 \begin{equation}
 r_\phi:=\frac 1{\sqrt{2\sinh \phi}}\begin{pmatrix}
 1&1\\e^{-\phi}&e^\phi
 \end{pmatrix}\in \SL_2(\R)
 .
 \label{eq:r-phi}
 \end{equation}

\begin{lemma}
The operator
\begin{multline}
J_\alpha^{-1}T_\alpha(r_\phi)J_\alpha f(t)
= (2\sinh \phi)^{\alpha/2} e^{-\phi(\alpha/2+it)}e^{\pi t/2}
\times\\ \times
\int_{-\infty}^\infty
f(s)\,
{}_2F_1 \left[\begin{matrix} \alpha/2-is,\alpha/2+it\\ \alpha\end{matrix};1-e^{-2\phi}  \right]
e^{-\pi s/2}
\mu(s)\,ds
.
\label{eq:JTJ}
\end{multline} 
is a unitary operator 
$L^2(\R,\mu(s)\,ds)\to L^2(\R,\mu(s)\,ds)$.
\end{lemma}
 
 {\sc Proof.}
 The operator 
 $$
 J_\alpha^{-1} T_\alpha(r_\phi) J_\alpha
 $$
 is unitary by definition as a product of 3 unitary operators
 $$
 L^2(\R,\mu(s)\,ds) \to H_\alpha\to H_\alpha \to L^2(\R,\mu(s)\,ds)
 .
 $$
 We must find explicit formula for composition.
   Write $J_\alpha$ in the form
 \begin{equation}
 J_\alpha f(z)=e^{i\pi\alpha/4}2^{\alpha}
 \int_{-\infty}^\infty f(s) z^{-\alpha/2-is} e^{-\pi s/2} d\mu(s)
, 
\label{eq:Jf}
\end{equation}
we use $(e^{-i\pi/2})^{-is}=e^{-\pi s/2}$. 
In this formula we take the branch of 
 $z^{-\alpha/2-is}$ given by
\begin{equation} 
  z^{-\alpha/2-is}=e^{-(\alpha/2+is)\ln z}
  \label{eq:lnz}
,  \end{equation}
 where the logarithm is real on the semi-axis
$z>0$. 
Then the inversion formula is
\begin{equation}
J_\alpha^{-1} F(t)
=
\frac{e^{-i\pi\alpha/4} 2^{-\alpha}\Gamma(\alpha)e^{\pi t/2}}
{\Gamma(\alpha/2+it) \Gamma(\alpha/2-it)} 
\int_0^{\infty} F(z) z^{\alpha/2+it-1}ds
\label{eq:J-}
.
\end{equation} 

Recall that $r_\phi$ is given by (\ref{eq:r-phi}),
\begin{multline*}
T_\alpha(r_\phi)J_\alpha f(z)=
 e^{i\pi\alpha/4}2^{\alpha}(2\sinh \phi)^{\alpha/2}
 \times\\ \times
\int_{-\infty}^\infty  \left(\frac {e^\phi z+1}{e^{-\phi}z+1}\right)^{\alpha/2-is}
(e^{-\phi}z+1)^{-\alpha} f(s)e^{-\pi s/2}\,d\mu(s)
=
\end{multline*}
\begin{multline}
=\frac{e^{i\pi\alpha/4}2^{\alpha}(2\sinh \phi)^{\alpha/2}}
{2\pi\Gamma(\alpha)}
 \times\\ \times
 \int_{-\infty}^\infty  (e^\phi z+1)^{-\alpha/2-is}
(e^{-\phi}z+1)^{-\alpha/2+is} f(s)\,e^{-\pi s/2}|\Gamma(\alpha+is)|^2ds
\label{eq:TJ}
.
\end{multline}

Next, we apply the inverse transform $J_\alpha^{-1}$,
\begin{multline}
J_\alpha^{-1}T_\alpha(r_\phi)J_\alpha f(t)
= \frac{(2\sinh \phi)^{\alpha/2} e^{\pi t/2}}
{\Gamma(\alpha/2+it) \Gamma(\alpha/2-it)}
\times\\ \times
\int_0^\infty \int_{-\infty}^\infty
z^{\alpha/2+it-1} (e^\phi z+1)^{-\alpha/2-is}
(e^{-\phi}z+1)^{-\alpha/2+is} f(s)
\times\\ \times
e^{-\pi s/2}
|\Gamma(\alpha+is)|^2ds\,dz
.
\label{eq:j-inverse}
\end{multline}
We must evaluate the integral in $z$,
\begin{multline*}
\int_0^\infty z^{\alpha/2+it-1} (e^\phi z+1)^{-\alpha/2-is}
(e^{-\phi}z+1)^{-\alpha/2+is}\,dz=\\=
e^{-\phi(\alpha/2+it)}
\int_0^\infty 
u^{\alpha/2+it-1} (1+u)^{-\alpha/2-is}
(1+e^{-2\phi}u)^{-\alpha/2+is}\,du
=\\=
e^{-\phi(\alpha/2+it)}\cdot
\frac{\Gamma(\alpha/2+it)\Gamma(\alpha/2-it)}{\Gamma(\alpha)} 
\,\,
{}_2F_1 \left[\begin{matrix} \alpha/2-is,\alpha/2+it\\ \alpha\end{matrix};1-e^{-2\phi}  \right]
,
\end{multline*}
here we applied an integral representation 
of the Gauss hypergeometric function,
\begin{equation}
_2F_1
\left[\begin{matrix} a,b\\c\end{matrix}; 1-u \right]
=\frac{\Gamma(c)}{\Gamma(b)\Gamma(c-b)}
\int_0^\infty y^{b-1}(1+y)^{a-c}(1+y u)^{-a}dy
\label{eq:F21}
,
\end{equation}
see \cite{HTF1}, 2.12(5); this is valid for $|\arg u|<\pi$.

Thus, we get (\ref{eq:JTJ}).
\hfill $\square$

\sm

{\sc Proof of Theorem \ref{th:vilenkin-2}.} 
Finally, we change function by the rule
\begin{equation}
\label{eq::}
g(s)=e^{-\pi s/2} f(s)
.
\end{equation}
This is equivalent to passing to the space $L^2(\R,w(s)\,ds)$,
where $w(s)=e^{\pi s} d\mu(s)$.
 
\sm 
 
%%%%%%%%%%%%%%%%%%%%%%%%%%%%%%%%%%%%%%%%%%%%%%%%%%%%%%%%%%%%%%%%%

{\bf\punct  Calculations. The difference operator.} 
Now we  evaluate the image of the operator 
$$
f(s)\mapsto s f(s)
 %, \qquad M:L^2(\R,d\mu(s))\to L^2(\R,d\mu(s)),
$$
under  $J_\alpha^{-1}T_\alpha(r_\phi)J_\alpha$.
Differentiating (\ref{eq:TJ}) by parameter $z$, we get
that  $J_\alpha^{-1}T_\alpha(r_\phi)$ send
the operator
$$
f(s)\mapsto -2is\sinh\phi f(s)
$$
to
$$
D:=
(z^2+2z\cosh \phi+1) \frac  d{dz}+\alpha(z+\cosh\phi) 
.$$

Next, we evaluate the corresponding operator in $L^2(\R,w(s)\,ds)$.
First, set 
\begin{equation}
g(t)=\int_0^\infty F(z) z^{\alpha/2+it-1}dz, \qquad
h(t)=\frac {g(t)}{\Gamma(\alpha/2+is)\Gamma(\alpha/2-is)} 
\label{eq:gh}
\end{equation}
and evaluate
\begin{multline}
\int_0^\infty DF(z) z^{\alpha/2+it-1}dz
=\\=
\int_0^\infty 
F'(z)\left(z^{\alpha/2+it+1} + 2\cosh \phi\, z^{\alpha/2+it} +z^{\alpha/2+it-1}
\right)\,dz+\\+ \alpha
\int_0^\infty F(z)\left(z^{\alpha/2+it} +
\cosh z^{\alpha/2+it-1}\right)\,dz
\label{eq:DF}
.
\end{multline}
Next we {\it formally} integrate by parts and come
to
$$
(\alpha/2-it-1)g(t-i)-2\cosh \phi\,t g(t) +(-\alpha/2-it-1)g(t+i)
.
$$
For functions $h\in L^2(\R,w(t)\,dt)$ we get the transformation
$$
h(t)\mapsto
(\alpha/2+it) h(t-i)- 2i \cosh \phi \,t h(t)- (\alpha/2-it) h(t+i)
.$$

%%%%%%%%%%%%%%%%%%%%%%%%%%%%%%%%%%%%%%%%%%%%%%%%%%%%%%%%%%%%%%%%%%%%%%%%%%%%%%%%%%%
 
 {\bf\punct Self-adjointness. Proof of Theorem \ref{th:vilenkin-1}.} 
  Denote by $W_R$ the space of functions $f(s)$ holomorphic in the strip
\begin{equation}
|\Im s|<R
\label{eq:strip}
\end{equation}
 satisfying the condition:
 for any $A>0$ there is $C$ such that
$$ 
|f(s)|< C\cdot \exp (-A|\Re s|)
.
$$ 
The operator $f\mapsto sf(s)$ in $L^2(\R,d\mu(s))$ 
is essentially self-adjoint 
on $W_R$.

Theorem \ref{th:vilenkin-1} is a corollary of the following lemma.

\begin{lemma}
If $R$ is sufficiently large, then for any $f\in W_R$
we have $\frV_\alpha f\in \cH[w]$.
\end{lemma}

 {\sc Proof.} 
 Since $f(z)$ super-exponentially decreases,
  $J_\alpha f(z)$, see (\ref{eq:Jf})--(\ref{eq:lnz}), is 
a well-defined
 analytic function on the universal covering of $\C\setminus 0$.
In other words, we can assume in (\ref{eq:lnz}) that
 $-\infty<\arg z<+\infty$. 
Since $f$ is analytic in the strip, the Fourier transform 
of $f$  exponentially 
decreases, therefore the Mellin transform decreases as $O(|z|^R)$ 
 as $|z|\to 0$
and as $O(|z|^{-R})$  as $z\to \infty$ (see \cite{Tit},  Theorem 31),
both $O(\cdot)$ are uniform in any  sector $|\arg z|<C$ with
 finite central angle.

 After the transform $T_\alpha(r_\phi)$ we get a function 
 $F(z):=T_\alpha(r_\phi)J_\alpha f(z)$ on the universal covering over
 $$
 \C\setminus \{-e^{-\phi},-e^{\phi}\}
.
 $$
It has the following behavior near the ramification points:

\sm

1. Near $\infty$ the function $F(z)$ 
 has form $z^{-\alpha}\gamma(1/z)$, where $\gamma$ is
holomorphic near $0$.

\sm

2. Near $e^\phi$ we have $F(z)=O(z-e^\phi)^R$.

\sm
 
3. Near $e^{-\phi}$ we have  $F(z)=O(z-e^{-\phi})^{R-\alpha}$.

\sm

Dominants $O(\cdot)$ are uniform in all sectors 
with finite central angles.
 
Next, we   examine the function $g(t)$ given by (\ref{eq:gh}).
The function $F(z)z^{\alpha/2}$ is holomorphic in the sector 
$|\arg z|<\pi$ and admit estimates $O(|z|^{\alpha/2})$ at zero 
and $O(|z|^{-\alpha/2})$ 
at $\infty$. Therefore (see \cite{Tit}, Theorem 31),
 its Mellin transform $g(t)$ is 
 
--- holomorphic in the strip
$|\Im t|<\alpha/2$, 

--- decreases as $O(e^{-(\pi-\epsilon) |\Re t|})$ as 
$\Re t\to\pm\infty$. 

Both consequences are not sufficient for our purposes%
\footnote{If $\alpha\le 2$, then the width of the strip 
is not sufficient. }.
For this reason, we improve a behavior of $F(z)$ at zero and at infinity
(in the spirit of Watson's Lemma%
\footnote{See, e.g., \cite{AAR}, Theorem C.3.1}).
 
 Consider the functions 
 \begin{align*}\tau_1(z)&=
 \exp(-z^{1/3})\left(1+z^{1/3}+\frac 1{2!}z^{2/3}\right),
 \\
 \tau_2(z)&=z^{-\alpha}
 \exp(-z^{-1/3})\left(1+z^{-1/3}+\frac 1{2!}z^{-2/3}\right)
 .
 \end{align*}
 
 \begin{lemma}
 The functions
$$ R(t):= \int_0^\infty \tau_1(z) z^{\alpha/2+it-1}dz,\qquad
Q(t):= \int_0^\infty \tau_2(z) z^{\alpha/2+it-1}dz
$$
are meromorphic in the strip 
$$
-\alpha/2-1<\Im t<\alpha/2+1.
$$
A unique singularity of $R(t)$ in the strip is a simple pole at $t=i\alpha$.
A unique singularity of $Q(t)$ in the strip is a simple pole at $t=-i\alpha$.
Both functions admit the following estimate in the strip
\begin{equation}
 O( |t|^{3\alpha/2-1/2} e^{-3\pi |t|/2}), \qquad \Re t\to\pm\infty
.
\label{eq:corr}
\end{equation}
 \end{lemma}

{\sc Proof.} 
$$
R(t)=3\left(\Gamma(3\alpha/2+3it)+ \Gamma(3\alpha/2+1+3it)+\frac 12\Gamma(3\alpha/2+2+3it)
\right)
. $$
 Poles of summands are $i\alpha/2$, $i\alpha/2+i/3$, $i\alpha/2+2i/3$, but the last
 two poles cancel.
 \hfill $\square$

\sm

Next, consider the function
$$
F^\circ(z):=F(z)-F(0)\tau_1(z)-\bigl(z^\alpha F(z)\bigr)\Bigl|_{z=\infty}
 \cdot \tau_2(z)
.
$$
Denote
$$
g(t)=\int_0^\infty z^{\alpha/2} F(z) z^{is-1}ds,
\qquad
g^\circ(t)=\int_0^\infty z^{\alpha/2} F^\circ(z) z^{is-1}ds,
$$

The  function $z^{\alpha/2}F^\circ(z)$  admits the following expansions
near $0$ and $\infty$
\begin{align}
z^{\alpha/2}F^\circ(z)&=
p_1z^{\alpha/2+1}+ p_2 z^{\alpha/2+2}+p_3 z^{\alpha/2+3} +\dots,
 \quad |z|\to 0,
\label{eq:add1}
\\
  z^{\alpha/2} F^\circ(z)&=
q_1z^{-\alpha/2-1}+ q_2 z^{-\alpha/2-2}+q_3 z^{-\alpha/2-3}
 +\dots,
,\quad  |z|\to\infty
\label{eq:add2}
\end{align}
in the sector $|\arg z|\le\pi$. It is  continuous up to the boundary
of the sector if $R>\alpha$.
The functions 
$$
\gamma_\pm(x):=z^{\alpha/2}F^\circ(z)\Bigr|_{z=e^{x\pm i\pi}}
$$
 have $R-\alpha$ derivatives. Expansions
(\ref{eq:add1})--(\ref{eq:add2}) imply the following lemma.
 
 \begin{lemma}
 \label{l:decreasing}
All derivatives  $\frac{d^k}{dx^k}\gamma_\pm(x)$ tend to zero as  $x\to\pm\infty$. 
\end{lemma}

Therefore (see \cite{Tit}, Theorem 31 and proof of Theorem 26),
 $g^\circ(t)$ is holomorphic in the strip
$|\Im t|<\alpha/2+1$ and satisfy the estimate
$$
|g^\circ(t)|=O(e^{-\pi |\Re t|} |\Re t|^{-(R-\alpha)}),
\qquad \Re t \to\pm\infty
.
$$

The function $g(t)$ satisfy the same estimate  (because $g(t)-g^\circ(t)$
is (\ref{eq:corr})) at infinity, but it is meromorphic in the strip with
simple poles at $t=\pm\alpha/2$.

Now it remains to divide%
\footnote{The formula (\ref{eq:j-inverse}) contains also a multiplication by
$e^{\pi t/2}$, but this factor cancels after (\ref{eq::}).}
 $g(t)$ by $\Gamma(\alpha/2+it)\Gamma(\alpha/2-it)$.
Poles at $t=\pm i\alpha/2$ disappear, we get a function holomorphic in the strip
$|\Im t|<1+\alpha/2$, and decreasing as $O(t^{-(R-2\alpha+1)})$.
It remains to choose a sufficiently wide strip (\ref{eq:strip}).
\hfill $\square$

%%%%%%%%%%%%%%%%%%%%%%%%%%%%%%%%%%%%%%%%%%%%%%%%%%%%%%%%%%%%%%%%%%%%%%%%%%%%%%%%
%%%%%%%%%%%%%%%%%%%%%%%%%%%%%%%%%%%%%%%%%%%%%%%%%%%%%%%%%%%%%%%%%%%%%%%%%%%%%%%%%
%%%%%%%%%%%%%%%%%%%%%%%%%%%%%%%%%%%%%%%%%%%%%%%%%%%%%%%%%%%%%%%%%%%%%%%%%%%%%%

\section{Example of self-adjoint extensions}

\COUNTERS

This section  contains another construction in Vilenkin's style, see
 \cite{VK}, Section 7.7.11. A representation-theoretic standpoint
 of our considerations is explained at the end of the section.
 
 %%%%%%%%%%%%%%%%%%%%%%%%%%%%%%%%%%%%%%%%%%%%%%%
 
 \sm

{\bf \punct The  difference operator.}
Consider the space $L^2(\R)$ and the subspace
 $\cV$ consisting of functions
holomorphic in the strip $|\Im s|\le 1$ and decreasing as
$$
|f(s)|=O(\Re s)^{-3/2-\epsilon},\qquad |\Re s|\to\infty
.
$$

Fix $\tau\in\R$, $0<\phi<\pi$ and consider the operator
$$
\cL f(s)=i(1/2-is) f(s+i)+ 2(s-\tau)\cos\phi f(s) -i(1/2+is-2it) f(s-i)
.
$$

It is symmetric, see Subsection 2.4.
 
\begin{proposition}
\label{l:1-1}
The  operator $\cL$ is not self-adjoint. Its defect indices are $(1,1)$. 
\end{proposition} 

The author does not know are self-adjoint extensions of $\cL$
 natural objects or not.
For this reason we consider another example. 
 
Consider the operator $\cL\oplus\cL$ acting in the space
$L^2(\R,ds)\oplus L^2(\R,e^{2\pi s}ds)$.

 Consider the space $\cH$ consisting of pair of functions 
$(f_1,f_2)$ meromorphic  in the strip
$|\Im s|\le 1$ such that 
\begin{align}
f_1(s)&= O(\Re s)^{-3/2-\epsilon},
\qquad |\Re s|\to\infty
\label{eq:f1f2-1}
\\
 e^{2\pi s} f_2(s)&= O(\Re s)^{-3/2-\epsilon},
\qquad |\Re s|\to\infty
\label{eq:f1f2-2}
.
\end{align}

Fix $\sigma\in\R$. Consider the space $\cH_\sigma$ consisting of pair of functions 
$(f_1,f_2)$ meromorphic  in the strip
$|\Im s|\le 1$ and satisfying (\ref{eq:}), with simple poles at points 
$i/2$ and $-i/2+2\tau$. We also require
\begin{align}
\res\limits_{s=i/2} f_1(s)&=
\res\limits_{s=i/2} f_2(s),
\label{eq:res1}
\\
\res\limits_{s=-i/2+2\tau} f_1(s)&=-e^{2\pi(\tau+i\sigma)}
 \res\limits_{s=-i/2+2\tau} f_2(s)
\label{eq:res2}
.
\end{align}
The parameter $\sigma$ is present only in the last condition, it is a parameter
of a self-adjoint extension.

\begin{proposition}
{\rm a)} The operator $\cL\oplus \cL$  has defect indices $(2,2)$ on $\cH$.

\sm

{\rm b)} The operator $\cL\oplus \cL$ is essentially self-adjoint 
 on the domain 
$\cH_\alpha$.
\label{l:ss-extension}
\end{proposition}

Next, consider the following elements of the space
$\cH_\alpha$:
\begin{equation} 
\left(\Psi_1^{(n)}(s),\Psi_2^{(n)}(s)\right),
\label{eq:Psi-Psi}
\end{equation} 
 where both functions $\Psi_1^{(n)}$, $\Psi_2^{(n)}$ are given by the same formula
 $$
 B(1/2+is,1/2+2i\tau-is)
\,\,
 {}_2F_1\left[\begin{matrix}1/2+is,1/2-\sigma+i\tau-n\\1+2i\tau \end{matrix}
; 1-e^{-2i\phi} 
 \right]
. $$
The function $\Psi_1^{(n)}$ is obtained by analytic continuation of 
$$
B(\dots)\,\,
 {}_2F_1\left[\begin{matrix}1/2+is,1/2-\sigma+i\tau-n\\1+2i\tau \end{matrix}
; z
 \right]
$$
 from $z=0$ along the path $z=1-e^{-2i\theta}$ with $\theta\in [0,\phi]$;
$\Psi_2^{(n)}$  along the path $z =1-e^{2i\theta}$ with $\theta\in [0,\pi-\phi]$.

\begin{proposition}
\label{pr:basis}
{\rm a) (\cite{Ner-mei})} Elements 
$\left(\Psi_1^{(n)}(s),\Psi_2^{(n)}(s)\right)$,
where $n$ ranges in $\Z$, form an orthogonal basis
in the space 
$L^2(\R,ds)\oplus L^2(\R,e^{2\pi s}ds)$.

\sm

{\rm b)} They also are eigenfunctions of the operator
 $\cL\oplus \cL$ defined on $\cH_\alpha$. The eigenvalues
 are $2\sin\phi(\sigma+n)$.
\end{proposition}

%%%%%%%%%%%%%%%%%%%%%%%%%%%%%%%%%%%%%%%%%%%%%%%%%%%%%%%%%%%%%%%%%%

\sm

{\bf\punct A family of orthogonal bases in $L^2(\R)$.} 
Fix $\tau\in\R$, $\sigma\in\C$ and $\phi\in (0,\pi)$.
Define functions
$$\Delta_{\sigma}(x)=\Delta_{\sigma}(x;\tau,\phi)
=
(1+xe^{i\phi})^{-1/2-i\tau-\sigma}
(1+xe^{-i\phi})^{-1/2-i\tau+\sigma}
.
$$
We choose a branch of $\Delta_{\sigma}(x)$ 
 by the condition $\Delta_{\sigma}(0)=1$.

\begin{lemma}
For any $\tau$, $\sigma\in\R$, the functions
 $\Delta_{\sigma+n}$,
where $n$ ranges in $\Z$, form an orthogonal basis in $L^2(\R)$.
\end{lemma}

{\sc Proof.} We pass to a new variable $\theta\in [0,2\pi]$
defined by
$$
e^{i\theta}=\frac{1+e^{i\phi}x}{1+e^{-i\phi}x},
\qquad d\theta=
\frac{2\sin\phi\, dx}{(1+e^{i\phi}x)(1+e^{-i\phi}x)}
.$$
Then
$$
(2\sin\phi)^{1/2+i\tau}\cdot \Delta_{\sigma+n}=
e^{-i(\sigma+n)\theta}
 \theta'(x)^{1/2+i\tau} 
$$
We consider the map from $L^2[0,2\pi]$ to
$L^2(\R)$ given by
$$
S f(x)=f(\theta(x))\theta'(x)^{1/2+i\tau}
$$
Evidently, it is unitary.
The system $\Delta_{\sigma+n}$ is the image of the
complete orthogonal
system $e^{-i(\sigma+n)\theta}$ under the map $S$.
\hfill $\square$

\sm

%%%%%%%%%%%%%%%%%%%%%%%%%%%%%%%%%%%%%%%%%%%%%%%%%%%%%%%%%%%%%%%%%%

{\bf \punct A differential operator.}
Fix $\tau\in\R$, $\phi\in(0,\pi)$. We consider
the following symmetric differential operator
\begin{equation}
D=D_{\tau,\phi}
=i(x^2+2\cos\phi x+1)\frac {d}{dx}+i(1+2i\tau)(x+\cos\phi)
\label{eq:Dtauphi}
\end{equation}
in $L^2(\R,dx/2\pi)$.

The functions $\Delta_{\sigma}$
  are formal eigenfunctions of the operator $D$,
\begin{equation}
D \Delta_{\sigma}(x)= (2\sin\phi)\,\sigma
\Delta_{\sigma}(x)
.
\end{equation}

\begin{lemma}
{\rm a)} Defect indices of the operator $D$ defined on the subspace
$C^\infty_c(\R)$ are $(1,1)$.

\sm

{\rm b)} Defect indices of the operator $D$ defined on the subspace
$C^\infty_c\bigl((0,\infty) \bigr)$ are $(1,1)$.
\end{lemma}

{\sc Proof.} Indeed, functions $\Delta_{\sigma}$ are contained in $L^2(\R)$
for all $\sigma\in\C$. Therefore, $\dim\ker(D^*\pm i)=1$.
\hfill $\square$

Fix $\sigma\in\R$.
Denote by $W_\alpha$ the space of $C^\infty$-functions on $\R$ such that there is
a function $h(y)$ smooth near zero such that
\begin{equation}
f(x)=\begin{cases}
x^{-1-2i\tau} h(1/x),& \text{for sufficiently large positive $x$}
\\
e^{-2\pi i\sigma}(-x)^{-1-2i\tau} h(1/x),& \text{for sufficiently small negative $x$}
\end{cases}
\label{eq:cases-h}
\end{equation}

\begin{lemma}
The operator $D$ is essentially self-adjoint on the subspace $W_\sigma$
and $\Delta_{\sigma+n}$ are its eigenfunctions.
\end{lemma}

{\sc Proof.} Verification of symmetry of $D$ on $W_\sigma$
 is straightforward.
 The subspace $W_\sigma$ contains vectors $\Delta_{\sigma+n}$.
Other functions $\Delta_{\kappa}$ are not in the domain of definiteness
of $D^*$ and therefore defect indices are $(0,0)$.
\hfill$\square$

%%%%%%%%%%%%%%%%%%%%%%%%%%%%%%%%%%%%%%%%%%%%%%%%%%%%%%%%%%%%%%%%%%

\sm

{\bf\punct The double  Mellin transform.} 
Let $f\in L^2(\R)$. Consider the pair of functions
\begin{align}
g_1(s)&=\int_0^\infty f(x) x^{is-1/2} dx
\label{eq:g1}
,
\\
g_2^\circ(s)&=\int_{-\infty}^0 f(x)(-x)^{is-1/2} dx
\label{eq:g2-circ}
.
\end{align}
Obviously,
$$
\int_{-\infty}^\infty |f(x)|^2 dx=
\frac 1{2\pi}
\left\{
\int_{-\infty}^\infty |g_1(s)|^2ds+
\int_{-\infty}^\infty |g_2^\circ(s)|^2ds
\right\}
.
$$

Thus we get a unitary operator 
$L^2(\R,dx)\to L^2(\R, ds/2\pi)\oplus L^2(\R, ds/2\pi)$.
Let modify this transform and set
\begin{equation}
g_2(s)=\int_{-\infty}^\infty f(x) x^{is-1/2}ds=-i
e^{-\pi s}\int_0^\infty f(x) (-x)^{is-1/2}ds
,
\label{eq:g2}
\end{equation}
here we take a branch of $x^{is-1/2}$ that is analytic in
 the upper half-plane
and real for $x>0$. Now we get
$$
\int_{-\infty}^\infty |f(x)|^2 dx=
\frac 1{2\pi}
\left\{
\int_{-\infty}^\infty |g_1(s)|^2ds+
\int_{-\infty}^\infty |g_2(s)|^2 e^{2\pi s}ds
\right\}
.
$$
We denote the operator $f\mapsto (g_1,g_2)$ by $\wt\frM$

%%%%%%%%%%%%%%%%%%%%%%%%%%%%%%%%%%%%%%%%%%%%%%%%%

\sm

{\bf\punct The difference operator.} We evaluate the $\wt\frM$-image 
of $Df$ as in (\ref{eq:DF}) and get the formal difference  operator
$\cL\oplus \cL$ in $L^2(\R,ds)\oplus L^2(\R,e^{2\pi s} ds)$.

 Propositions \ref{l:1-1},
 \ref{l:ss-extension}.a are corollaries of the following lemma.

\begin{lemma}
{\rm a)} The image of $C^\infty_{c}(0,\infty)$ under (\ref{eq:g1}) is contained in
$\cV$.

\sm

{\rm b)} The image of $C^\infty_{c}(-\infty,0)+C^\infty_{c}(0,\infty))$
is contained in $\cH$.
\end{lemma}

{\sc Proof.} a) Recall that the Mellin transform of $f$ is reduced
to the Fourier transform by the substitution
$x=e^y$ to $f(x)$. In (\ref{eq:g1}) 
we evaluate the Fourier transform of
$f(e^y)e^{y/2}$,
the function $g_1(s)$ decreases as $O(s^{-N})$ for any $N$.
  
  \sm  
  
b) We apply the same argument to $g_2^\circ(s)$, see (\ref{eq:g2-circ}).
 After passing to $g_2$ we get the estimate
 (\ref{eq:f1f2-2}).
 \hfill $\square$

\sm

Proposition \ref{l:ss-extension}.b is a corollary of the following lemma.

\begin{lemma}
The image of the space $W_\alpha$ under the Mellin transform
$\wt\frM$ is contained in the space $\cH_\alpha$.
\end{lemma}  

{\sc Proof.} We  repeat considerations in the spirit 
of Watson lemma.
Pass to the function
$$
f^\star(x)=
\begin{cases}
f(x)- f(0)e^{-x}- h(0)x^{-1-2i\tau}e^{-1/x},
 \qquad & x>0
\\
f(x)-f(0)e^{x}-h(0) e^{-2\pi i\sigma} (-x)^{-1-2i\tau} e^{1/x}, 
\qquad & x<0
\end{cases}
.
$$
where $h$ is the same as in (\ref{eq:cases-h}).
Consider the first component
of the transform $\wt\frM$.
 
 We have
\begin{align}
f^\star(x)&=  c_1 x+\dots + c_N x^N+ O(x^{N+1}), \qquad
 &x\to 0+
\\
f^\star(x)&= d_1 x^{-2-2i\tau} +\dots + d_M x^{-M-2i\tau}+ O(x^{-M-1}),
\qquad & x\to +\infty
.\end{align}

Examine the behavior of
$$
g_1(s)=\int_0^\infty f(x) x^{is-1/2} dx,
\qquad
g^\star_1(s)=\int_0^\infty f^\star(x) x^{is-1/2} dx 
.
$$
Functions $g_1(s)$, $g_1^\star(s)$
are Fourier transforms of $f(e^y)e^{y/2}$, $f^\star(e^y)e^{y/2}$.
It is easy to see that derivatives of $f^\star(e^y)e^{y/2}$
admit estimates
$$
\frac{d^k}{dy^k} \bigl( f^\star(e^y)e^{y/2} \bigr)=
O(e^{-3|y|/2})
.
$$
Therefore $g_1^\star(s)$ is defined in the strip
$|\Im s|<3/2$ and decreases in this strip as $O(|\Re s|^{-N})$
for any $N$.

On the other hand, 
$$g_1(s)-g_1^\star(s)=f(0)\Gamma(1/2+is)+h(0)\Gamma(-1/2-2i\tau+is)
$$
is meromorphic in the stir with poles at
$s=i/2$, $s=-i/2+2i\tau$
and
exponentially decreases as $|\Re s|\to\infty$.
The residues at poles are $f(0)$ and $h(0)$ respectively.

In the same way we prove decreasing of $g_2^\circ(s)$ at infinity.
Residues at poles $s=i/2$, $s=-i/2+2i\tau$ are respectively $f(0)$
and $e^{-2\pi\sigma i}$. It remains to multiply 
$g_2^\circ(s)$ by $-ie^{-\pi s} h(0)$ and we come
to (\ref{eq:res1})--(\ref{eq:res2}).
\hfill $\square$

\sm

%%%%%%%%%%%%%%%%%%%%%%%%%%%%%%%%%%%%%%%%%%%%%%%

{\bf\punct Proof of Proposition \ref{pr:basis}.}
We evaluate $\wt\frM \Delta_{\sigma+n}$
using the formula (\ref{eq:F21}) and come to (\ref{eq:Psi-Psi}).

%%%%%%%%%%%%%%%%%%%%%%%%%%%%%%%%%%%%%%%%%%%%%%%%%%%%%%%

\sm
 
{\bf\punct The origin of construction of this section.}
Fix $\sigma$, $\tau\in\R$.
Consider the following representation
$T_{\tau,\sigma}(g)$ 
 of
the group $\SL_2(\R)$ in $L^2(\R)$,
$$
T_{\tau,\sigma}\begin{pmatrix}a&b\\c&d \end{pmatrix}
f(x)=
f\left(\frac{b+xd}{a+zc}  \right)
(a+zc)^{-1/2-\sigma+i\tau} \ov{\ln(a+zc)}^{\,-1/2+\sigma+i\tau} 
$$
In this formula, we choose any branch of $\ln(a+zc)$
that is holomorphic in the upper half-plane
and define powers as
\begin{multline*}
(a+zc)^{-1/2-\sigma+i\tau} \ov{(a+zc)}^{\,-1/2+\sigma+i\tau} 
:=\\:=
\exp\left((-1/2-\sigma+i\tau)\ln(a+zc)+ (-1/2+\sigma+i\tau)  \ov{(a+zc)}\right)
\end{multline*}
Thus, an operator $T_{\tau,\sigma}(g)$
 is determined up to a constant factor
and we get a projective unitary representation of
$\SL(2,\R)$ (it is a representation of the principal series,
 see, e.g.,
\cite{Ner-gauss},  Subsection 7.4.3).

The operator $D_{\tau,\phi}$ given by (\ref{eq:Dtauphi})
is an infinitesimal generator of the group $\SL_2(\R)$. It generates
a compact subgroup, and $\Delta_{\sigma+n}$ are eigenvectors of this subgroup.

The transform $\wt\frM$ is the
 spectral decomposition of the one-parametric group of operators
 $T_{\tau,\sigma}\begin{pmatrix} a&0\\0&a^{-1}  \end{pmatrix}$.

%%%%%%%%%%%%%%%%%%%%%%%%%%%%%%%%%%%%%%%%%%%%%%%%%%%%%%%%%%%%%%%%%%%%

%<vibogach@mail.ru>

{\tt Math.Dept., University of Vienna,

 Nordbergstrasse, 15,
Vienna, Austria

\&

Institute for Theoretical and Experimental Physics,

Bolshaya Cheremushkinskaya, 25, Moscow 117259,
Russia

\&

Mech.Math. Dept., Moscow State University,
Vorob'evy Gory, Moscow

e-mail: neretin(at) mccme.ru

URL:www.mat.univie.ac.at/$\sim$neretin

wwwth.itep.ru/$\sim$neretin}


\small

\begin{thebibliography}{cc}

\bibitem{AAR}
Andrews, G. E.; Askey, R.; Roy, R. {\it Special functions.}
 Cambridge University Press, Cambridge, 1999.

\bibitem{AW0}
Askey, R.; Wilson, J. {\it A set of hypergeometric orthogonal polynomials.} 
 SIAM J. Math. Anal.  13  (1982), no. 4, 651--655.

\bibitem{AW}
Askey, R.; Wilson, J. {\it Some basic hypergeometric orthogonal polynomials
 that generalize Jacobi polynomials.} 
 Mem. Amer. Math. Soc.  54  (1985),  no. 319.

\bibitem{Che1}
Cherednik, I.
{\it Inverse Harish-Chandra transform and difference operators.}
Internat. Math. Res. Notices 1997, no. 15, 733--750. 

\bibitem{Che2}
Cherednik, I. {\it Double affine Hecke algebras.}
 Cambridge University Press, Cambridge, 2005.

\bibitem{Car}
Carlitz, L. {\it Bernoulli and Euler numbers and orthogonal polynomials.}
  Duke Math. J  26  (1959) 1--15

\bibitem{Dunford}
Dunford, N.; Schwartz, J. T. {\it Linear operators. Part II: Spectral theory.
 Self adjoint operators in Hilbert space.}
 John Wiley \& Sons, New York-London, 1963 

\bibitem{HTF1}
Erd\'elyi, A.; Magnus, W.; Oberhettinger, F.; Tricomi, F. G.
{\it Higher transcendental functions.} Vol. I.
Based, in part, on notes left by Harry Bateman. McGraw-Hill Book Company,
 Inc., New York-Toronto-London, 1953.



\bibitem{HTF2}
Erd\'elyi, A.; Magnus, W.; Oberhettinger, F.; Tricomi, F. G.
{\it Higher transcendental functions.} Vol.  II.
Based, in part, on notes left by Harry Bateman. McGraw-Hill Book Company,
 Inc., New York-Toronto-London, 1953

%\bibitem{HCh}

\bibitem{Gro}
Groenevelt, W.
{\it The Wilson function transform.}
Int. Math. Res. Not. 2003, no. 52, 2779--2817. 

\bibitem{Her}
H\"ormander, L. {\it The analysis of linear partial
 differential operators. I. Distribution theory 
and Fourier analysis.}
 Springer-Verlag, Berlin, 1983.



\bibitem{Koe2}
Koekoek, R.; Lesky, P. A.; Swarttouw, R. F.
{\it Hypergeometric orthogonal polynomials and their $q$-analogues.}
 Springer-Verlag, Berlin, 2010.

\bibitem{Koe}
Koekoek, R.;
Swarttouw, R. F.
{\it The Askey-scheme of hypergeometric orthogonal polynomials
 and its $q$-analogue.}
Delft University of Technology
Faculty of Information Technology and Systems
Department of Technical Mathematics and Informatics.
Report no. 98-17,
1998.  Available via {\tt http://aw.twi.tudelft.nl/$\sim$koekoek/askey/}


\bibitem{Koo}
Koornwinder, T.
{\it A new proof of a Paley-Wiener type theorem for the Jacobi transform.}
Ark. Mat. 13 (1975), 145--159. 


\bibitem{KoLe}
 Lebedev, N.N.,  Kontorovich, M.I.,
{\it On the application of inversion formulae to the solution
 of some electrodynamics problems,} 
J. Exper. Theor. Phys. 9(6) (1939), pp. 729--742 (in Russian). 

\bibitem{Leb}
Lebedev, N. N.; Skalskaya, I. P.; Uflyand, Y. S. 
{\it Problems of mathematical physics.}
Prentice-Hall, Inc.,
 Englewood Cliffs, N.J., 1965;
Reprinted as {\it Worked problems in applied mathematics,}
Dover, 1979


\bibitem{Les1}
Lesky, P. A. {\it Unendliche und endliche Orthogonalsysteme 
von continuous Hahnpolynomen.}
 (German)   Results Math.  31  (1997),  no. 1-2, 127--135. 


\bibitem{Les2}
Lesky, P. A.; Waibel, B.
{\it  Orthogonalit\"at von Racahpolynomen und Wilsonpolynomen.}
 (German) 
  Results Math.  35  (1999),  no. 1-2, 119--133.

\bibitem{Mei}
Meixner, J. {\it Orthogonale Polynomsysteme mit einer besonderen
 Gestalt der erzeugenden Funktion.} (German)
 J. of  London Math. Soc. 9, 1934, 6--13. 

\bibitem{Ner-ber}
Neretin, Yu. A.
{\it The index hypergeometric transform and an imitation 
of the analysis of Berezin kernels on hyperbolic spaces.}
 Sb. Math.  192  (2001),  no. 3-4, 403--432.




\bibitem{Ner-wilson}
Neretin, Yu. A. 
{\it Beta integrals and finite orthogonal systems of Wilson polynomials.}
 Sb. Math.  193  (2002),  no. 7-8, 1071--1089. 

\bibitem{Ner-mei}
Neretin, Yu. A. {\it 
Perturbations of classical hypergeometric orthogonal systems.}
Addendum to preprint
{\tt arXiv:math/0309445} (2003).

\bibitem{Ner-double}
Neretin, Yu. A. {\it Some continuous analogues of the expansion
 in Jacobi polynomials, and vector-valued orthogonal bases.}
  Funct. Anal. Appl.  39  (2005),  no. 2, 106--119.



\bibitem{Ner-gauss}
Neretin, Yu. A. {\it 
Lectures on Gaussian integral operators and classical groups.}
European. Math. Soc, 2011.

\bibitem{NUS}
Nikiforov, A. F.; Suslov, S. K.; Uvarov, V. B.
{\it Classical orthogonal polynomials of a discrete variable.}
 (Russian) 
Nauka, Moscow, 1985. 216 pp. There is an extended English version,
Springer, 1991.

\bibitem{Ole}
Olevski\u i, M. N. 
{\it On the representation of an arbitrary function
 in the form of an integral with a kernel containing 
a hypergeometric
 function.} 
(Russian) Doklady Akad. Nauk SSSR (N.S.)  69,  (1949). 11--14. 


\bibitem{PBM1}
Prudnikov, A. P.; Brychkov, Yu. A.; Marichev, O. I. 
{\it Integrals and series. Vol. 1. Elementary functions.}
 Gordon and Breach, New York, 1986.

\bibitem{PBM3}
Prudnikov, A. P.; Brychkov, Yu. A.; Marichev, O. I. 
{\it Integrals and series. Vol. 3. More special functions.}
 Gordon and Breach,
 New York, 1990. 

\bibitem{Rom}
 Romanovski, V. I. {\it Sur quelques classes nouwels 
of polynomes orthogonaux.} (French)
Compt. Rend. Acad. Sci. Paris, 188 (1929), 1023--1025.

\bibitem{Sla}
Slater, L. J. {\it Confluent hypergeometric functions.}
 Cambridge University Press, New York, 1960

\bibitem{Tit}
Titchmarsh, E. C.
{\it Introduction to the theory of Fourier integrals.}
Third edition. Chelsea Publishing Co., New York, 1986. 


\bibitem{Vil}
Vilenkin, N. Ja. 
{\it Special functions and the theory of group representations.}
American Mathematical Society, Providence, R. I., 1968.

\bibitem{VK}
Vilenkin, N. Ja.; Klimyk, A. U. 
{\it Representation of Lie groups and special functions. Vol. 1.
 Simplest Lie groups, special functions and integral transforms.}
 Kluwer, Dordrecht, 1991

\bibitem{Wat}
Watson, G. N. {\it A Treatise on the Theory of Bessel Functions.}
 Cambridge University Press, New York, 1944.

\bibitem{Wil}
Wilson, J. A.
{\it Some hypergeometric orthogonal polynomials.}
SIAM J. Math. Anal. 11 (1980), no. 4, 690--701.

\bibitem{Wim}
Wimp, J. {\it A class of integral transforms.} 
 Proc. Edinburgh Math. Soc. (2),  14  (1964/1965) 33--40.

\bibitem{YaLu}  
 Yakubovich, S. B. {\it Index transforms.}
World Scientific, 1996


\end{thebibliography}
\end{document}